\documentclass[VANCOUVER,STIX1COL]{WileyNJD-v2}

\usepackage{float}
\usepackage{epsfig, graphicx}
\usepackage{latexsym,amsfonts,amsbsy,amssymb}
\usepackage{amsmath,amscd}
\usepackage[capitalize]{cleveref}
\usepackage{bm}
\usepackage{tikz, pgfplots}
\usetikzlibrary{arrows}
\usepackage{ulem}
\usepackage{placeins}

\renewcommand{\vec}[1]{\mathbf{#1}}

\articletype{Article Type}%

\received{26 April 2016}
\revised{6 June 2016}
\accepted{6 June 2016}

\begin{document}

\title{A positivity-preserving unigrid method for elliptic PDEs}

\author[1]{Ronald D. Haynes}

\author[1]{Scott MacLachlan}

\author[2]{Dawei Wang}

\authormark{Haynes, MacLachlan and Wang}

\address[1]{\orgdiv{Department of Mathematics and Statistics}, \orgname{Memorial
  University of Newfoundland}, \orgaddress{\state{Newfoundland}, \country{Canada}}}

\address[2]{\orgdiv{Department of Computer Science}, \orgname{University of Toronto}, \orgaddress{\state{Ontario}, \country{Canada}}}


\corres{*\email{smaclachlan@mun.ca}}


\abstract[Summary]{
  While constraints arise naturally in many physical models, their treatment in
  mathematical and numerical models varies widely, depending on the nature of
  the constraint and the availability of simulation tools to enforce it.  In
  this paper, we consider the solution of discretized PDE models that have a
  natural constraint on the positivity (or non-negativity) of the solution.
  While discretizations of such models often offer analogous positivity
  properties on their exact solutions, the use of approximate solution
  algorithms (and the unavoidable effects of floating--point arithmetic) often destroy any guarantees that the
  computed approximate solution will satisfy the (discretized form of the)
  physical constraints, unless the discrete model is solved to much higher
  precision than discretization error would dictate.  Here, we introduce a class
  of iterative solution algorithms, based on the unigrid variant of multigrid
  methods, where such positivity constraints can be preserved throughout the
  approximate solution process.  Numerical results for one- and two-dimensional
  model problems show both the effectiveness of the approach and the trade-off required to ensure positivity of approximate solutions throughout the solution process.
}



\maketitle

\section{Introduction}\label{sec:intro}


Many real-world problems in science and engineering have ``constrained''
solution sets, where the physical model is only sensible under certain
conditions.  The treatment of such constraints in computational science and
engineering varies dramatically depending on the nature of the constraint.  Some
constraints are preserved strongly in the mathematical modeling and
discretization processes, some are preserved only in a weak sense, and some are
altogether discarded, either for computational convenience or for more
principled reasons.  In this paper, we consider constraints on positivity (or
non-negativity) of solutions to PDEs, which arise in many settings.  Positive or
non-negative solutions arise naturally in Markov Chain
models\cite{schweitzer1991survey,stewart1992numerical,krieger1995numerical,de2010algebraic,de2010smoothed,treister2010},
diffusion of mass or chemical concentration\cite{jovanovic2008}, and thickness
of thin films\cite{bertozzi1998,matta2011}.  They also arise indirectly in some
cases, such as in adaptive mesh refinement in 1D, where a non-tangling
constraint on the mesh leads to a positivity constraint on the
meshwidths\cite{white1979}.  In many of these models, the physical constraint on
positivity is reflected in the discrete model, such as through stochastic or
M-matrix properties\cite{berman1994nonnegative}, guaranteeing that the exact
solution of the discrete model should satisfy a discrete version of the
constraint.

Such guarantees are reassuring but, in practice, difficult to rely upon due to
the inherent nature of floating-point arithmetic and iterative solution
algorithms.  Quite simply, we \textit{never} solve discrete systems exactly and,
so, properties of the exact solution are never guaranteed in our discrete
approximations.   In this paper, we consider the question of how to guarantee
positivity of inexact solutions using multigrid algorithms, a well-known class of efficient algorithms for the solution of many discretized elliptic PDEs\cite{briggs2000multigrid, UTrottenberg_etal_2001a}.  Standard (geometric) multigrid methods achieve fast convergence by progressively eliminating oscillatory error components within a hierarchical decomposition of the problem and, hence, the coarse-grid correction process can easily introduce negative components into an approximate solution.  The main contribution herein is recognizing that pointwise positivity constraints can naturally be enforced in such a hierarchical setting, but not within a classical V-cycle algorithm.

To our knowledge, the only application area where efficient positivity-preserving algorithms are in widespread use is in that of computing the stationary probability distribution of a Markov Chain.  The family of iterative aggregation-disaggregation (IAD) techniques  was first proposed in the 1990s and has been extensively studied\cite{schweitzer1986iterative,schweitzer1991survey,stewart1992numerical, krieger1995two}. Following this work, efficient multigrid methods\cite{krieger1995numerical,de2010smoothed,de2010algebraic,treister2010} were proposed. One of the key factors in the success of these algorithms is that the system matrices involved are singular M-matrices and, hence, a regular splitting can be used as a relaxation algorithm that preserves positivity of approximate solutions. In addition, the prolongation/restriction operators correspond to aggregation/disaggregation operators that have meaningful probabilistic interpretations.  Thus, while these techniques have been very successful for models leading to Markov Chains, they cannot be readily adapted to more general problems, such as finite-difference or finite-element discretizations of PDEs, due to the lack of appropriate prolongation/restriction operators.

Within the standard ``correction scheme'' approach to multigrid, enforcing
fine-level positivity constraints on coarse-level computations is difficult,
since we must allow negative corrections to fine-level approximations that
overshoot the true solution.  To overcome this roadblock, we propose to make use
of the unigrid method introduced by McCormick and Ruge
\cite{mccormick1983unigrid, WHolland_SFMcCormick_JWRuge_1982a}, to understand
the effects of fine-level constraints on the coarser levels of the hierarchy.
Unigrid is an algorithmic framework that results from interpreting the steps in
a multigrid algorithm by their direct effect on the finest-grid approximation.
That is, each change from the coarse-grid correction procedure is immediately
propagated to the finest-grid approximation and, thus, used in subsequent steps
of the coarse-grid correction.  Under certain assumptions, unigrid and classical
multigrid are equivalent, although they have different computational
characteristics.  In practice, unigrid algorithms are more expensive than their
multigrid counterparts, but they can offer advantages for implementation and convergence.
 Unigrid methods can also be derived
from other viewpoints, such as that of subspace-correction methods
\cite{griebel2012greedy,manteuffel2006projection}.  This was used in the
projection multilevel method of Manteuffel et al.~\cite{manteuffel2006projection}, where subspace
correction was used to avoid direct linearization in the solution of nonlinear
problems.
This viewpoint has also been used to present randomized and greedy analogues of the multigrid process\cite{griebel2012greedy,MGriebel_POswald_2020a}.

The main contribution of this paper is to recognize that the unigrid framework provides a natural avenue to directly enforce fine-level constraints on the coarse-grid correction process within a multigrid algorithm.  In this way, we derive and present an approach where the finest-grid approximate solution can be guaranteed to remain positive throughout the solution process, subject to natural requirements on the relaxation scheme on all grids.  Moreover, by combining the robustness of algebraic multigrid methods with this unigrid strategy, we are able to develop a very efficient and generally-applicable iterative approach that can preserve positivity of approximate solutions. 

The remainder of the paper is organized as follows.  In \cref{sec:model_problems}, we first present our model problems and discuss related questions regarding preserving positivity of their approximation solutions. In \cref{sec:unigrid}, we introduce the central ideas of the unigrid methodology, and formulate the main algorithms of our positivity-preserving unigrid methods. To test the efficiency of these algorithms, we give numerical results on the model problems in \cref{sec:numerical_experiments}, followed by concluding remarks in \cref{sec:conclusion}.

\section{Model problems and discretization}\label{sec:model_problems}

While many physical problems may include positivity constraints, in this paper we consider only those that arise from relatively simple elliptic model problems in one and two spatial dimensions.  For these problems, standard discretization approaches yield discrete systems whose exact solutions naturally reflect the positivity at the discrete level.

\subsection{Model Problems}\label{ssec:model}

We consider the solution of the one- and two-dimensional diffusion equations
\begin{subequations}
\begin{align}
  -\frac{d}{dx}\left( a(x,u)\frac{du}{dx} \right) + b(x)u & = 0 \text{ for }0<x<1, \label{eq:1d_nonlinear} \\
  -\nabla \cdot \left(a(x,y) \nabla u\right) & = f(x,y)\text{ for }(x,y)\in(0,1)^2,
\end{align}
\end{subequations}
with prescribed Dirichlet boundary conditions.  We assume that the diffusion coefficient, $a$, is positive-valued and bounded below by a constant strictly greater than zero.  Under suitable conditions on the zeroth-order term, $b$, and the Dirichlet boundary data, positivity of the continuum solution to such an equation can be guaranteed using, for example, a maximum principle argument \cite{LCEvans_2010a}.

As a first example, one of the applications that inspired the work of this paper
is the equidistributing mesh generation problem in
1D\cite{WHuang_RDRussell_2011a}. A mesh transformation $u(x)$ which
equidistributes a prescribed mesh density function, $a(u)$, is the solution of
the nonlinear two--point boundary value problem
\begin{equation*}
    \frac{d}{dx}\bigg( a(u)\frac{du}{dx} \bigg) = 0,\;u(0)=\theta_0,\;u(1)=\theta_1,
\end{equation*}
where $0$ and $1$ are the boundaries of the computational domain and $\theta_0$ and $\theta_1$ are the boundaries of the physical domain.  Direct calculation shows that the exact (continuum) solution of this problem is monotonically increasing, leading to untangled meshes on the one-dimensional interval $[\theta_0,\theta_1]$, and thus, the mesh spacings are positive.  Preliminary experiments showed that standard numerical methods for discretizing this equation and solving the resulting linear system did not preserve this property.

\subsection{Positivity-preserving discretization}
\label{ssec:discretization}

As is well-recognized in many settings, properties of the continuum solution of PDEs need not be preserved in the discretization process.  Here, we identify some characteristics of discretizations that will be important for the success of our positivity-preserving algorithm.  Corresponding to nonnegative/positive functions that are useful in the continuous case, we first give a definition of a nonnegative/positive matrix (or, as a special case with $n=1$, a vector). 
\begin{definition}
Let $A = [a_{ij}]$ be a real $m\times n$ matrix.  We say $A \geq 0\; (A > 0)$ if $a_{ij} \geq 0\; (a_{ij} > 0)$ for all $1 \leq i \leq m,\; 1 \leq j \leq n$, and call $A$ a nonnegative (positive) matrix. Similarly, $A = 0$ means that every entry of $A$ is 0.
\end{definition}

Suppose the discretization of a linear continuum PDE gives a linear system $A\mathbf u = \mathbf b$, with $A$ being nonsingular. Since the solution to this system is $\mathbf u = A^{-1}\mathbf b$, an obvious sufficient condition to have a nonnegative solution is to require $\mathbf b \geq 0$ and $A^{-1} \geq 0$. As has been long-established, M-matrices have a close connection with this property\cite{berman1994nonnegative}. For convenience of discussion, we adopt the common notation and let $Z^{n\times n}$ be the set of square matrices with non-positive off-diagonal and nonnegative diagonal entries, i.e. matrices of the form
\[
A = 
\begin{bmatrix}\label{matrix:Z_matrix}
a_{11}  & -a_{12} & -a_{13} & \cdots & -a_{1n}\\
-a_{21} & a_{22}  & -a_{23} & \cdots & -a_{2n} \\
-a_{31} & -a_{32} & a_{33}  & \cdots & -a_{3n} \\
\vdots  & \vdots  & \vdots  & \ddots & \vdots \\
-a_{n1} & -a_{n2} & -a_{n3} & \cdots & a_{nn}
\end{bmatrix},
\]
where $a_{ij} \geq 0$. 
An important subclass of $Z$-matrices are the M-matrices, defined as follows.
\begin{definition}
Any $n\times n$ real matrix $A$ that can be written in the form $A = sI-B$ with $s>\rho(B)$ and $B\geq0$ is called an M-matrix, where $\rho(B)$ is the spectral radius of $B$. 
\end{definition}
One important property of M-matrices is that if $A$ is an M-matrix, then $A$ is nonsingular and $A^{-1}\geq 0$.  This can be strengthened to $A^{-1}>0$ if $A$ is also irreducible~\cite{berman1994nonnegative}.  In this case, if \textit{any} entry in $\mathbf b$ is positive, then $\mathbf u = A^{-1} \mathbf b$ is a (strictly) positive vector.

Considering the PDE models of the previous section, it is clear that in the case of $b(x) = 0$, standard finite-difference or finite-element discretizations of the linear model problems lead to M-matrices. If $b(x)$ is a non-negative function, then a finite-difference discretization of the PDE also leads to an irreducible M-matrix structure, but the situation is more complicated for finite-element discretizations (due to positive off-diagonal entries in the standard mass matrix).  For the nonlinear case, Newton linearization of the model problems leads to the inclusion of first-order derivative terms on $u$, which are difficult to discretize in a manner that preserves M-matrix structure.  To avoid this complication, we consider simple Picard linearization techniques that yield slower convergence to the solution of the nonlinear problem, but allow us to use positive approximate solutions of the preceding linearizations to ensure positivity of the exact solution to the linearized system at each step of the calculation. To see this, consider the centered finite-difference of the derivative term in~\cref{eq:1d_nonlinear} on a uniform mesh 
\begin{equation*}
    -\frac{d}{dx}\bigg( a(x,u)\frac{du}{dx} \bigg)_j
    \approx \frac{1}{h^2} \bigg( -a(x_{j-1/2},u_{j-1/2})u_{j-1} + \big( a(x_{j-1/2},u_{j-1/2}) + a(x_{j+1/2},u_{j+1/2}) \big) u_j - a(x_{j+1/2},u_{j+1/2}) u_{j+1} \bigg),
\end{equation*}
where $u_{j} = u(x_{j})$ and so on. The coefficient values can be approximated by 
\begin{equation*}
    a(x_{j-1/2},u_{j-1/2}) = a\big(x_{j-1/2}, (u_{j-1}+u_j)/2 \big), \quad a(x_{j+1/2},u_{j+1/2}) = a\big(x_{j+1/2}, (u_{j}+u_{j+1})/2 \big).
\end{equation*}
Then, discretizing the rest of~\cref{eq:1d_nonlinear} in a consistent manner yields a nonlinear algebraic system of the form $A(\mathbf u) \mathbf u - \mathbf f(\mathbf u) = \mathbf 0$.
Several possible approaches can be used to linearize this system, including a Picard iteration or Newton's method. To use Newton's method, we write the functional
\begin{equation*}
    \mathbf F(\mathbf u) = A(\mathbf u) \mathbf u - \mathbf f(\mathbf u) = \mathbf 0.
\end{equation*}
Then, for a given initial guess $\mathbf u^{(0)}$, Newton's method derives subsequent approximations to $\mathbf u$ by computing
\begin{equation*}
    \mathbf u^{(k+1)} = \mathbf u^{(k)} + \mathbf s^{(k)},\; \text{where}\; \mathbf F'(\mathbf u^{(k)})\mathbf s^{(k)} = \mathbf F(\mathbf u^{(k)}),
\end{equation*}
for $k\geq 0$, where $ \mathbf F'(\mathbf u^{(k)})$ is the Jacobian matrix.  Of particular interest here is that $\mathbf F'(\mathbf u^{(k)})$ not only depends on the system matrix, $A(\mathbf u)$, but also the right-hand side $\mathbf f(\mathbf u)$ and their derivatives (with respect to $\mathbf u$). Therefore, it is not guaranteed that $F'(\mathbf u^{(k)})$ is a $Z$-matrix, nor that $\mathbf F(\mathbf u^{(k)})$ is a non-negative vector.  Indeed, for an arbitrary positive approximation, $\mathbf u^{(k)}$, it is not even clear that positivity of $\mathbf s^{(k)}$ is a desirable property, since it may be that $\mathbf u^{(k)} \geq \mathbf u^{(k+1)} > 0$.

Thus, to have a chance at a positivity-preserving method, it seems more feasible to use a Picard iteration for linearization, which means to replace the unknown $\mathbf u$ in $A(\mathbf u)$ and $f(\mathbf u)$ with values from the previous iteration, i.e.,
\begin{equation} \label{eqn:picard}
    A(\mathbf u^{(k-1)})\mathbf u^{(k)} = \mathbf f(\mathbf u^{(k-1)}) = \mathbf b^{(k-1)}.
\end{equation}
Note that  (when the problem is discretized as above) $A(\mathbf u^{(k-1)})$ is a nonsingular irreducible M-matrix, so the iteration is well-posed and the exact solution is positive when the nonzero right-hand side satisfies $\mathbf b^{(k-1)} \geq 0$ (with at least one strictly positive entry). 
\cref{alg:Picard_iteration} gives the pseudocode for a Picard iteration, in which the solution step inside the while loop is the main topic of discussion in this paper. We will call the Picard iteration the outer loop, and the loop to solve for the next iterative approximation $\mathbf u_1$ the inner loop.

\begin{algorithm}
  \begin{algorithmic}[1]
    \Procedure{PICARD}{$A,\mathbf u,\mathbf f$}
\While{$||A(\mathbf u)\mathbf u-\mathbf f(\mathbf u)|| > \tau$}
\State Solve $A(\mathbf u)\mathbf u_1 = \mathbf f(\mathbf u)$ for  $\mathbf u_1$
\State Set $\mathbf u = \mathbf u_1$
\EndWhile
\State Return $\mathbf u$
\EndProcedure
\caption{Picard Iteration}
\label{alg:Picard_iteration}
\end{algorithmic}
\end{algorithm}


\section{Multigrid and Unigrid}\label{sec:unigrid}

Consider the solution of the linear system $A\vec{u} = \vec{b}$.  From the
multigrid perspective, we let $A^{(0)} = A$ be the ``finest-grid'' matrix, and
assume that we are given a hierarchy of discrete operators and grid-transfer
operators, that we denote by $A^{(k)}$ for $0 \leq k \leq l$ and $I_k^{k+1}$
(for restriction from grid $k$ to grid $k+1$) and $I_{k+1}^k$ (for interpolation
from grid $k+1$ to grid $k$).  In geometric multigrid, the grid-transfer
operators can be defined, for example, as canonical finite-element interpolation
operators or their (scaled) counterparts in the finite-difference case, and the
coarse-grid operators can be determined by rediscretization.  For algebraic
multigrid, the interpolation and restriction operators on each level are
determined by the discrete operator on that level, and the Galerkin product,
$A^{(k+1)} = I_k^{k+1}A^{(k)}I_{k+1}^k$, is used to define the operator on the
next level.  With these operators, we can define a multigrid V-cycle as in
\cref{alg:Vcycle}, where each level of the hierarchy stores a single vector,
$\vec{u}^{(k)}$, as the approximate solution to the residual-correction equation
on that level.  The multigrid algorithm naturally relies on two smoothing or
relaxation operators, $M_1^{(k)}$ and $M_2^{(k)}$, used for pre- and
post-relaxation, respectively.  Standard choices include the inverse of the (weighted) diagonal of $A^{(k)}$, leading to the (weighted) Jacobi iteration, or of the lower-triangular part of $A^{(k)}$, leading to the Gauss-Seidel iteration.  Parameters $\nu_1$ and $\nu_2$ denote the number of pre- and post-relaxation sweeps to be performed on each level.  For comparison to the unigrid algorithm presented next, we consider the case where $\nu_2 = 0$ and, thus, drop the subscript on $\nu_1$ in what follows.

\begin{algorithm}
  \begin{algorithmic}[1]
    \Procedure{VCYCLE}{$\mathbf u^{(0)},\mathbf b^{(0)}, \nu_1, \nu_2, \{A^{(k)}\}, \{M_1^{(k)}\}, \{M_2^{(k)}\}, \{I_{k}^{k+1}\}, \{I_{k+1}^k\}$}
  \For{$k = 0, 1, \ldots, l-1$}
  \State Do $\nu_1$ steps of the iteration, $\vec{u}^{(k)} = \vec{u}^{(k)} +
    M_1^{(k)}\left(\vec{b}^{(k)}-A^{(k)}\vec{u}^{(k)}\right)$, with initial guess
    $\vec{u}^{(k)}$
  \State Compute $\vec{b}^{k+1} = I_{k}^{k+1}\left(\vec{b}^{(k)}-A^{(k)}\vec{u}^{(k)}\right)$, set $\vec{u}^{(k+1)} =
    \vec{0}^{(k+1)}$
    \EndFor
  \State Solve $A^{(l)}\vec{u}^{(l)} = \vec{b}^{(l)}$
  \For{$k = l-1, l-2, \ldots 0$}
  \State Do $\nu_2$ steps of the iteration, $\vec{u}^{(k)} = \vec{u}^{(k)} +
    M_2^{(k)}\left(\vec{b}^{(k)}-A^{(k)}\vec{u}^{(k)}\right)$, with initial guess
    $\vec{u}^{(k)}+I_{k+1}^k\vec{u}^{(k+1)}$
    \EndFor
\State Return $\mathbf u^{(0)}$
\EndProcedure
  \caption{Multigrid V-cycle}\label{alg:Vcycle}
  \end{algorithmic}
\end{algorithm}

\subsection{Unigrid background}

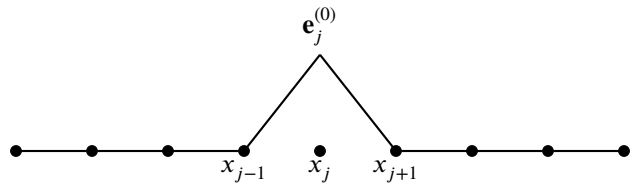
\begin{figure}
\begin{center}
\begin{tikzpicture}
    \filldraw (-4 cm, 0) circle (2pt) node[align=center, above] {};
    \filldraw (-3 cm, 0) circle (2pt) node[align=center, above] {};
    \filldraw (-2 cm, 0) circle (2pt) node[align=center, above] {};
    \filldraw (-1 cm, 0) circle (2pt) node[align=center, below] {$x_{j-1}$};
    \filldraw (0 cm, 0) circle (2pt) node[align=center, below] {$x_j$};
    \filldraw (1 cm, 0) circle (2pt) node[align=center, below] {$x_{j+1}$};
    \filldraw (2 cm, 0) circle (2pt) node[align=center, above] {};
    \filldraw (3 cm, 0) circle (2pt) node[align=center, above] {};
    \filldraw (4 cm, 0) circle (2pt) node[align=center, above] {};
    
    \draw [thick] (-4 cm, 0) -- (-1 cm, 0) node[midway]{};
    \draw [thick] (0 cm, 36pt) -- (-1 cm, 0) node[midway]{};
    \draw [thick] (0 cm, 36pt) -- (1 cm, 0) node[midway]{};
    \draw [thick] (1 cm, 0) -- (4 cm, 0) node[midway]{};
    \filldraw (0 cm, 36pt) circle (0pt) node[align=right, above] {$\mathbf e_j^{(0)}$};
\end{tikzpicture}
\end{center}
\caption{Basis function, $\vec{e}_j^{(0)}$, associated with finest-mesh discretization.} \label{fig:GS_dir}
\end{figure}

A natural approach to unigrid comes from the perspective of a ``directional iteration'' \cite{mccormick1983unigrid,manteuffel2006projection} interpretation of the standard Gauss-Seidel iteration. Consider the $j$-th step of Gauss-Seidel iteration for solving $A\vec{u} = \vec{b}$, which corrects the $j$-th component of current approximate solution by an amount $\delta_j^{(0)}$
\begin{equation*}
    \mathbf u \leftarrow \mathbf u + \delta_j^{(0)}\mathbf e_j^{(0)},
\end{equation*}
where $\mathbf e_j^{(0)}$ is the $j$-th column of the $n\times n$ identity matrix, and
\begin{equation*}
    \delta_j^{(0)} = \frac{\langle\mathbf b-A\mathbf u, \mathbf e_j^{(0)}\rangle}{\langle A\mathbf e_j^{(0)}, \mathbf e_j^{(0)}\rangle},
\end{equation*}
such that the $j$-th component of residual is zero after the correction. The Gauss-Seidel method repeats this process for $1\leq j \leq n$, correcting one component of the approximate solution vector at each step until convergence. It is easy to see that the propagation of corrections in the Gauss-Seidel method is slow because of the "narrow" shape of the direction vectors $\mathbf e_j^{(0)}$, as shown in \cref{fig:GS_dir}, where we use $x_j$ to denote the $j^{\text{th}}$ point on a 1-dimensional discretization mesh. 

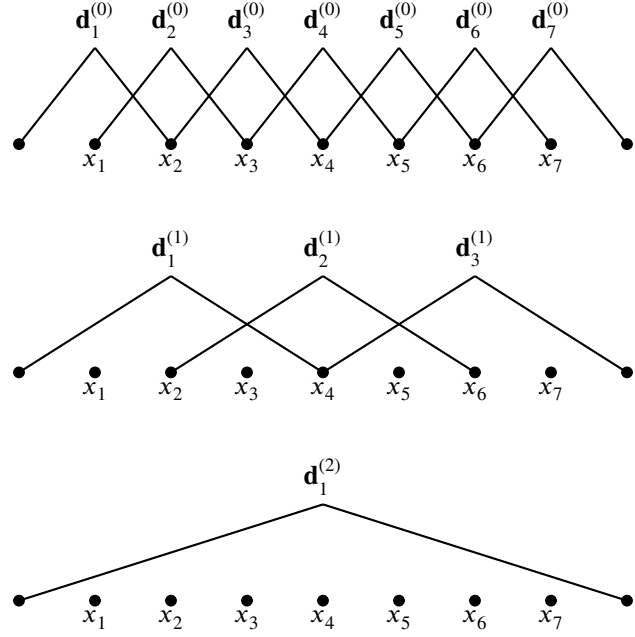
\begin{figure}
\begin{center}
\begin{tikzpicture}
    \filldraw (-4 cm, 0) circle (2pt) node[align=center, above] {};
    \filldraw (-3 cm, 0) circle (2pt) node[align=center, below] {$x_1$};
    \filldraw (-2 cm, 0) circle (2pt) node[align=center, below] {$x_2$};
    \filldraw (-1 cm, 0) circle (2pt) node[align=center, below] {$x_{3}$};
    \filldraw (0 cm, 0) circle (2pt) node[align=center, below] {$x_4$};
    \filldraw (1 cm, 0) circle (2pt) node[align=center, below] {$x_{5}$};
    \filldraw (2 cm, 0) circle (2pt) node[align=center, below] {$x_6$};
    \filldraw (3 cm, 0) circle (2pt) node[align=center, below] {$x_7$};
    \filldraw (4 cm, 0) circle (2pt) node[align=center, above] {};

    \draw [thick] (-3 cm, 36pt) -- (-4 cm, 0) node[midway]{};
    \draw [thick] (-3 cm, 36pt) -- (-2 cm, 0) node[midway]{};
    \filldraw (-3 cm, 36pt) circle (0pt) node[align=right, above] {$\mathbf d_1^{(0)}$};
    
    \draw [thick] (-2 cm, 36pt) -- (-3 cm, 0) node[midway]{};
    \draw [thick] (-2 cm, 36pt) -- (-1 cm, 0) node[midway]{};
    \filldraw (-2 cm, 36pt) circle (0pt) node[align=right, above] {$\mathbf d_2^{(0)}$};
    
    \draw [thick] (-1 cm, 36pt) -- (-2 cm, 0) node[midway]{};
    \draw [thick] (-1 cm, 36pt) -- (0 cm, 0) node[midway]{};
    \filldraw (-1 cm, 36pt) circle (0pt) node[align=right, above] {$\mathbf d_3^{(0)}$};
    
    \draw [thick] (0 cm, 36pt) -- (-1 cm, 0) node[midway]{};
    \draw [thick] (0 cm, 36pt) -- (1 cm, 0) node[midway]{};
    \filldraw (0 cm, 36pt) circle (0pt) node[align=right, above] {$\mathbf d_4^{(0)}$};
    
    \draw [thick] (1 cm, 36pt) -- (0 cm, 0) node[midway]{};
    \draw [thick] (1 cm, 36pt) -- (2 cm, 0) node[midway]{};
    \filldraw (1 cm, 36pt) circle (0pt) node[align=right, above] {$\mathbf d_5^{(0)}$};
    
    \draw [thick] (2 cm, 36pt) -- (1 cm, 0) node[midway]{};
    \draw [thick] (2 cm, 36pt) -- (3 cm, 0) node[midway]{};
    \filldraw (2 cm, 36pt) circle (0pt) node[align=right, above] {$\mathbf d_6^{(0)}$};
    
    \draw [thick] (3 cm, 36pt) -- (2 cm, 0) node[midway]{};
    \draw [thick] (3 cm, 36pt) -- (4 cm, 0) node[midway]{};
    \filldraw (3 cm, 36pt) circle (0pt) node[align=right, above] {$\mathbf d_7^{(0)}$};
    
    \filldraw (-4 cm, -3cm) circle (2pt) node[align=center, above] {};
    \filldraw (-3 cm, -3cm) circle (2pt) node[align=center, below] {$x_1$};
    \filldraw (-2 cm, -3cm) circle (2pt) node[align=center, below] {$x_2$};
    \filldraw (-1 cm, -3cm) circle (2pt) node[align=center, below] {$x_{3}$};
    \filldraw (0 cm, -3cm) circle (2pt) node[align=center, below] {$x_4$};
    \filldraw (1 cm, -3cm) circle (2pt) node[align=center, below] {$x_{5}$};
    \filldraw (2 cm, -3cm) circle (2pt) node[align=center, below] {$x_6$};
    \filldraw (3 cm, -3cm) circle (2pt) node[align=center, below] {$x_7$};
    \filldraw (4 cm, -3cm) circle (2pt) node[align=center, above] {};

    \draw [thick] (-2 cm, -3cm+36pt) -- (-4 cm, -3cm) node[midway]{};
    \draw [thick] (-2 cm, -3cm+36pt) -- (0 cm, -3cm) node[midway]{};
    \filldraw (-2 cm, -3cm+36pt) circle (0pt) node[align=right, above] {$\mathbf d_1^{(1)}$};
    
    \draw [thick] (0 cm, -3cm+36pt) -- (-2 cm, -3cm) node[midway]{};
    \draw [thick] (0 cm, -3cm+36pt) -- (2 cm, -3cm) node[midway]{};
    \filldraw (0 cm, -3cm+36pt) circle (0pt) node[align=right, above] {$\mathbf d_2^{(1)}$};
    
    \draw [thick] (2 cm, -3cm+36pt) -- (0 cm, -3cm) node[midway]{};
    \draw [thick] (2 cm, -3cm+36pt) -- (4 cm, -3cm) node[midway]{};
    \filldraw (2 cm, -3cm+36pt) circle (0pt) node[align=right, above] {$\mathbf d_3^{(1)}$};
    
    \filldraw (-4 cm, -6cm) circle (2pt) node[align=center, above] {};
    \filldraw (-3 cm, -6cm) circle (2pt) node[align=center, below] {$x_1$};
    \filldraw (-2 cm, -6cm) circle (2pt) node[align=center, below] {$x_2$};
    \filldraw (-1 cm, -6cm) circle (2pt) node[align=center, below] {$x_{3}$};
    \filldraw (0 cm, -6cm) circle (2pt) node[align=center, below] {$x_4$};
    \filldraw (1 cm, -6cm) circle (2pt) node[align=center, below] {$x_{5}$};
    \filldraw (2 cm, -6cm) circle (2pt) node[align=center, below] {$x_6$};
    \filldraw (3 cm, -6cm) circle (2pt) node[align=center, below] {$x_7$};
    \filldraw (4 cm, -6cm) circle (2pt) node[align=center, above] {};
    
    \draw [thick] (0 cm, -6cm+36pt) -- (-4 cm, -6cm) node[midway]{};
    \draw [thick] (0 cm, -6cm+36pt) -- (4 cm, -6cm) node[midway]{};
    \filldraw (0 cm, -6cm+36pt) circle (0pt) node[align=right, above] {$\mathbf d_1^{(2)}$};
\end{tikzpicture}
\end{center}
\caption{Three-level hierarchy of correction directions induced by linear interpolation from coarse grids.} \label{fig:unigrid_dir}
\end{figure}

Unigrid complements standard Gauss-Seidel relaxation by introducing ``broader-based'' directions, $\mathbf d^{(k)}_j$, defined on levels $1\leq k \leq l$. As illustrated in \cref{fig:unigrid_dir} for the case when there are three levels in 1D, a simple choice for the directions for levels $1\leq k \leq l$ are as the canonical unit vectors on those levels
 linearly interpolated to the finest grid, defined recursively by
\begin{equation}\label{eq:recursive_linear}
    \mathbf d_j^{(k)} = 
    \begin{cases}
    0.5\mathbf d^{(k-1)}_{2j-1} + \mathbf d_{2j}^{(k-1)} + 0.5\mathbf d^{(k-1)}_{2j+1}, & 1\leq k \leq l,\\
    \mathbf e_j^{(0)}, & k=0.
    \end{cases}
\end{equation}
Using these directions, one step of the unigrid update takes the form
\begin{equation}\label{eqn:unigrid_update}
    \mathbf u \leftarrow \mathbf u + \delta_j^{(k)} \mathbf d_j^{(k)},
\end{equation}
where now
\begin{equation*}
    \delta_j^{(k)} = \frac{\langle\mathbf b-A\mathbf u, \mathbf d^{(k)}_j\rangle}{\langle A\mathbf d_j^{(k)}, \mathbf d_j^{(k)}\rangle},
\end{equation*}
so that the resulting fine-grid residual is orthogonal to the direction $\mathbf d_j^{(k)}$ after correction. We can also add a weighting to the correction term, i.e., write
\begin{equation}\label{eqn:unigrid_step}
    \mathbf u \leftarrow \mathbf u + \omega\delta_j^{(k)}\mathbf d_j^{(k)},
\end{equation}
to damp the iteration.
A convenient way to denote the correction directions used within unigrid is to collect the vectors $\mathbf d_j^{(k)}$ as the columns of a matrix,
\begin{equation}\label{eqn:unigrid_directions}
    I_k^0 = [\mathbf d_1^{(k)}, \mathbf d_2^{(k)}, \cdots, \mathbf d_{n_k}^{(k)}], \; 0 \leq k \leq l,
\end{equation}
where $n_k$ is the total number of iteration directions on level $k$, and $l$ is the index for the coarsest level.
One way to view this is as the interpolation of the canonical unit vectors on level $k$ to the finest grid.  Note that the recursive definition used above in \cref{eq:recursive_linear} naturally leads to the relationship that $I_k^0 = I_1^0I_2^1\cdots I_k^{k-1}$ where each $I_k^{k-1}$ is just a standard multigrid interpolation operator (the 1D linear interpolation, in \cref{eq:recursive_linear}) from level $k$ to level $k-1$.

The relationship between directions for the unigrid iteration and interpolation operators in \cref{eqn:unigrid_directions} allows a more direct connection between a Gauss-Seidel iteration (over the directions within a level) and unigrid. Given any interpolation operators $I_1^0$, $I_2^1$, \ldots $I_{l}^{l-1}$, define $I_k^0 = I_1^0I_2^1\cdots I_{k}^{k-1}$ (with, trivially, $I_0^0$ being the identity operator on the finest grid).  The unigrid corrections associated with level $k$ of the hierarchy are given by 
\begin{align}
  {\delta}^{(k)}_j & = \frac{\langle \mathbf b-A\mathbf u, I_k^0\mathbf{e}_j^{(k)} \rangle}{\langle AI_k^0\mathbf{e}_j^{(k)}, I_k^0\mathbf{e}_j^{(k)}\rangle},\label{eqn:unigrid_delta} \\
  \vec{u} & \leftarrow \vec{u} + {\delta}^{(k)}_jI_k^0\mathbf{e}_j^{(k)},
\end{align}
for $1 \leq j \leq n_k$, the column dimension of $I_k^0$, where $\mathbf e_j^{(k)}$ is the $j$-th column of the $n_k\times n_k$ identity matrix.  With this, we can construct the unigrid analogue of any multigrid cycling strategy, simply by iterating over levels, $k$, as in the multigrid cycle.  Note also the strong connection between Gauss-Seidel on the restricted problem on level $k$ and the unigrid directional correction, best exposed by writing
\[
  {\delta}^{(k)}_j = \frac{\langle \mathbf b-A\mathbf u, I_k^0\mathbf{e}_j^{(k)} \rangle}{\langle AI_k^0\mathbf{e}_j^{(k)}, I_k^0\mathbf{e}_j^{(k)}\rangle} = \frac{\langle \left(I_k^0\right)^T\left(\mathbf b-A\mathbf u\right), \mathbf{e}_j^{(k)} \rangle}{\langle \left(I_k^0\right)^TAI_k^0\mathbf{e}_j^{(k)}, \mathbf{e}_j^{(k)}\rangle},
  \]
where we recognize the diagonal elements of the Galerkin operator on level $k$, $\left(I_k^0\right)^TAI_k^0$ in the denominator, and the restricted residual, $\left(I_k^0\right)^T\left(\mathbf b-A\mathbf u\right)$, in the numerator. Pseudocode for the analogue of a multigrid $V(\nu,0)$ cycle is given in \cref{alg:unigrid}.  Note that other fixed cycling strategies can be expressed in similar ways, or that greedy or randomized versions, that iterate over the columns of $[I,I_1^0,I_2^0,\ldots I_l^0]$, can also be defined \cite{griebel2012greedy, MGriebel_POswald_2020a}.

\begin{algorithm} 
  \begin{algorithmic}[1]
    \Procedure{UG\_SOLVE}{$A, \mathbf b, \mathbf u, \nu, \omega, \{I_k^0\}$} 
\For{$k = 0, 1, ..., l$} 
  \For{$i = 1, ..., \nu$} 
     \For{$j = 1, ..., n_k$} 
       \State Compute correction ${\delta}^{(k)}_j = \frac{\langle \mathbf b-A\mathbf u, I_k^0\mathbf{e}_j^{(k)} \rangle}{\langle AI_k^0\mathbf{e}_j^{(k)}, I_k^0\mathbf{e}_j^{(k)}\rangle}$
       \State Update approximation $\mathbf u = \mathbf u + \omega {\delta}^{(k)}_jI_k^0\mathbf{e}_j^{(k)}$
   \EndFor 
 \EndFor 
\EndFor
\State Return $\mathbf u$
\EndProcedure
\caption{Unigrid analogue of multigrid $V(\nu,0)$ cycle}\label{alg:unigrid}
\end{algorithmic}
\end{algorithm}
\vspace{5mm}

\subsection{Algebraic multigrid background}\label{sec:AMG}

The efficiency of multigrid methods lies in the fact that errors that are not effectively damped by relaxation on the finest grid can be reduced on coarser grids, through less-expensive relaxation on those grids.  While this is readily done for simple problems, where the space of slow-to-converge errors can be characterized by geometric smoothness, direct characterization of errors requiring coarse-grid correction is difficult in many cases.  Algebraic multigrid methods (AMG) offer an effective approach for some problems where such geometric characterization fails \cite{briggs2000multigrid}.  The key ``addition'' to geometric multigrid in AMG is a setup phase, in which the intergrid transfer operators and coarse-level operators are constructed recursively based on the given finest-grid matrix, defining the components of the multigrid cycle without any explicit dependence on information regarding the problem being solved other than the discretization matrix itself.  While many variations on AMG exist in the literature, we focus here on so-called ``classical'' (or ``Ruge-St\"uben'') AMG, as is commonly used as a solver for scalar PDEs like those considered here \cite{briggs2000multigrid, ABrandt_SFMcCormick_JWRuge_1984a, JWRuge_KStuben_1987a}.

Given the fine-grid matrix, $A$, the first step in the algebraic multigrid setup phase is determining a set of points to form the coarse grid.  Let $\Omega = \{1,2,\ldots,n\}$ denote the nodes of the mesh, used to label the (undirected) graph of connections in the matrix $A$, with edges, $(i,j)$, in the graph corresponding to nonzero entries, $a_{ij}\neq 0$, in the matrix.  AMG uses graph algorithms to determine $C\subset \Omega$ such that the points in $C$ make a ``good'' coarse-grid set for the matrix, $A$, where the heuristic definition of a good coarse-grid set is one from which an effective interpolation operator can be defined.  The determination of $C$ comes from a three-step process.  First, the graph of $A$ is filtered, removing edges corresponding to entries in $A$ that are small relative to the other entries in the same row.  Given a threshold parameter, $0 < \theta \leq 1$, an edge in the graph is said to be a strong connection if
\[
-a_{ij} \geq \theta \max_{k\neq i} \{-a_{ik}\},
\]
where we note the implicit assumption that $A$ is an M-matrix (or nearly so), so
that the off-diagonal entries are (predominantly) negative.  As is common in AMG~\cite{JWRuge_KStuben_1987a}, particularly for two-dimensional discretizations of diffusion-type operators, we take $\theta = 0.25$.  Any edge that does not satisfy this condition is then discarded as a ``weak connection'' in the matrix.  The next stage in the AMG coarsening algorithm is a first pass at determining the set $C$.  Here, a greedy algorithm is used to define a maximal independent set over the graph of strong connections.  Recall that an independent set of vertices in a graph is one where there is no edge in the graph between two vertices in the set, and that an independent set is maximal when no vertex can be added while preserving the property of being an independent set.  There is an intuitive nature to using a maximal independent set as the coarse set, $C$, since all points in $F = \Omega \setminus C$ possess at least one strong connection to a point in $C$.  However, there is no guarantee that a high-quality interpolation operator can be defined from any such set, $C$, since it may be the case that some points in $F$ have strong connections to other points in $F$ that cannot be easily accounted for in the interpolation process.  Thus, the final step in the AMG coarsening algorithm is to augment the preliminary $C$ set by examining each point in $F$ to ensure that any two strongly connected $F$ points have a \textit{common} (strongly-connected) neighbor in $C$, in order to ensure our ability to define a ``good'' interpolation operator.  In practice, this ``second pass'' adds a small number of points to the preliminary $C$ set in order to ensure this heuristic is satisfied.

Once the $C$ set is determined, the interpolation operator, $I_1^0$, is determined row-wise.  For rows, $i$, where $i \in C$, each row of $I_1^0$ has a single nonzero entry in it, with a value of 1 in the column corresponding to $i$'s index on the coarse mesh.  For rows $i\in F$, we partition the set $\{j\neq i | a_{ij}\neq 0\}$ into three disjoint subsets: $C_i$, the set of strongly-connected $C$ points, $F_i$, the set of strongly-connected $F$ points, and $W_i$, the set of weakly-connected neighbors.  The entries in row $i$ of $I_1^0$ are determined by first assuming that we seek to interpolate corrections to errors, $\vec{e}$, that satisfy $\left(A\vec{e}\right)_i \approx 0$, and then rewriting this as
\begin{equation}\label{eq:AMG_smallresidual}
a_{ii}e_i = -\sum_{j\in C_i}a_{ij}e_j - \sum_{k\in F_i}a_{ik}e_k - \sum_{\ell\in W_i}a_{i\ell}e_\ell.
\end{equation}
The first term on the right-hand side of this sum can be directly used in determining an interpolation operator, but we seek to eliminate the other two terms to determine an equation that relates just $e_i$ and $\{e_j\}_{j\in C_i}$.  To eliminate the sum over $F_i$, we make use of the assumption that each point $k\in F_i$ has a strong connection in $C_i$, ensured by the second pass of the $C$-selection algorithm.  Thus, we can approximate $e_k \approx \sum_{j\in C_i} a_{kj}e_j / \sum_{m\in C_i} a_{km}$, knowing that the sum in the denominator will not be zero (again assuming M-matrix-type properties) and that at least one value $a_{kj}$ will be significant.  For the final term, we make the approximation that $e_\ell \approx e_i$ for $\ell\in W_i$.  This, in fact, is not likely to be a good approximation; however, the connections in $W_i$ should not be significant, so any consistent approximation of $e_\ell$ should be acceptable.  Substituting these into \eqref{eq:AMG_smallresidual} yields
\[
\left(a_{ii} + \sum_{\ell\in W_i}a_{i\ell}\right)e_i = -\sum_{j\in C_i} \left(a_{ij} + \frac{\sum_{k\in F_i} a_{ik}a_{kj}}{\sum_{m\in C_i} a_{km}}\right)e_j.
\]
Dividing through by the coefficient on the left-hand side expresses $e_i$ as a linear combination of $\{e_j\}_{j\in C_i}$; these weights become the interpolation weights in row $i$ of $I_{1}^{0}$, in the columns corresponding to $j\in C_i$.

AMG then uses recursion and the Galerkin condition to define the remaining components in the multigrid cycle.  From $I_1^0$, we take $I_0^1 = \left(I_1^0\right)^T$, define $A^{(1)} = I_0^1 A^{(0)}I_1^0$, and then proceed recursively to define $I_2^1$ using $A^{(1)}$ in the role of the fine-grid matrix.  Because the graph-based coarsening algorithm does not reduce the grid size at a fixed rate, it is difficult to predict the number of levels in the multigrid hierarchy based on the fine-grid size.  In practice, we coarsen until the coarsest-grid matrix, $A^{(l)}$, is small enough that further coarsening would result in a coarsest grid of only 1 point.  In 1D, this typically leads to coarsest grids with 3 points (not including Dirichlet boundary conditions), but can lead to larger coarsest grids in 2D.  While it is typical in AMG to then directly factor this coarsest-grid matrix, we will consider only relaxation on the coarsest grid, as this is more compatible with the unigrid framework.

\subsection{Positivity-preserving unigrid}

In this section, we introduce two approaches that allow us to combine
unigrid with AMG-style interpolation in a way that guarantees preservation of
positivity of the approximate solution.  Broadly, these algorithms follow~\cref{alg:unigrid}, with the interpolation operators, $I_k^0$, being supplied by an AMG setup phase applied to the fine-grid discretization matrix.  Positivity of the approximate solution is maintained by post-processing the fine-grid correction, or the corrected fine-grid approximation, after each step of the unigrid algorithm.

\subsubsection{Uniform thresholding}
The main advantage of the unigrid approach in our context is the ability to check directly if a given correction, $\delta_j^{(k)}I_k^0\vec{e}_j^{(k)}$, will lead to negative entries in the fine-grid approximation after correction.  One approach to ensuring positivity is preserved in the approximate solution is, then, to weight corrections by always choosing a relaxation weight, $\omega$, such that
\begin{equation}\label{eq:unif_threshold_requirement}
    \mathbf u = \mathbf u + \omega{\delta}_j^{(k)}I_k^0\vec{e}_j^{(k)} > \mathbf 0,
\end{equation}
where the inequality is considered component-wise.  Because we use a multiplicative form of correction, we consider a ``default'' weight, $\omega = 1$, and only decrease $\omega$ if there are large-enough negative entries in the correction such that \eqref{eq:unif_threshold_requirement} would not be satisfied.  When the inequality in \eqref{eq:unif_threshold_requirement} would be violated with $\omega = 1$, we instead compute a weight guaranteed to ensure that it is satisfied with a smaller value for $\omega$.  To do this, we look over all entries of $\delta_j^{(k)}I_k^0\vec{e}_j^{(k)}$ that are negative, and compare these to the (positive) entries of $\vec{u}$ that they are correcting.  Defining $\mathcal{M}_i = \left\{m \middle| \left(\delta_j^{(k)}I_k^0\vec{e}_j^{(k)}\right)_m < 0\right\}$, enforcing \eqref{eq:unif_threshold_requirement} naturally requires that
\[
\omega < \min_{m\in \mathcal{M}_i} -u_m/\left(\delta_j^{(k)}I_k^0\vec{e}_j^{(k)}\right)_m.
\]
In practice, we take
\[
\omega = (1-\varepsilon) \min_{m\in \mathcal{M}_i} -u_m/\left(\delta_j^{(k)}I_k^0\vec{e}_j^{(k)}\right)_m,
\]
for a suitably small value of $\varepsilon$, such as $10^{-4}$.  This ensures
that the updated approximation stays positive, while still allowing for values
to be close to zero.  \Cref{alg:unigrid_uniform} gives the pseudocode of this variant of the unigrid algorithm.

\begin{algorithm} [h]
  \begin{algorithmic}[1]
    \Procedure{UG\_SOLVE\_THRESHOLDING}{$A, \mathbf b, \mathbf u, \nu,\varepsilon, \{I_k^0\}$}
\For{$k = 0, 1, ..., l$} 
  \For{$i = 1, ..., \nu$} 
     \For{$j = 1, ..., n_k$} 
       \State Compute correction ${\delta}^{(k)}_j = \frac{\langle \mathbf b-A\mathbf u, I_k^0\mathbf{e}_j^{(k)} \rangle}{\langle AI_k^0\mathbf{e}_j^{(k)}, I_k^0\mathbf{e}_j^{(k)}\rangle}$\\
       \If{$\min\{ \mathbf u + {\delta}^{(k)}_jI_k^0\mathbf{e}_j^{(k)} \} \leq 0$}
           \State Mark negative entries $\mathcal{M}_i = \left\{m \middle| \left(\delta_j^{(k)}I_k^0\vec{e}_j^{(k)}\right)_m < 0\right\}$
           \State Compute weight $\omega = (1-\varepsilon) \min_{m\in \mathcal{M}_i} -u_m/\left(\delta_j^{(k)}I_k^0\vec{e}_j^{(k)}\right)_m$
       \State Update approximation $\mathbf u = \mathbf u + \omega {\delta}^{(k)}_jI_k^0\mathbf{e}_j^{(k)}$
      \Else 
       \; Update approximation $\mathbf u = \mathbf u + {\delta}^{(k)}_jI_k^0\mathbf{e}_j^{(k)}$
       \EndIf
   \EndFor 
 \EndFor 
\EndFor
\State Return $\mathbf u$
\EndProcedure
\caption{Unigrid Method with Uniform Thresholding}\label{alg:unigrid_uniform}
\end{algorithmic}
\end{algorithm}

\subsubsection{Local correction}
Since the uniform thresholding technique relies on a single parameter, $\omega$, and scales all entries in the correction by the same amount (possibly forced by only a few points), it is possible that this might hamper the convergence because of a ``bad'' choice of $\omega$ that is constrained by some local behaviour. An alternative approach to avoid this possibility is to correct those negative components only locally.  In particular, because we expect that the solution varies smoothly, it makes sense to just locally correct any negative entries based on positive points ``nearby''.  Here, we propose two approaches, one based on interpolating values from positive-valued neighboring points, and one based on using local relaxation over negative entries until positivity is achieved.

A sketch of the first approach for a (hypothetical) one-dimensional problem is given in  \cref{fig:unigrid_local_correction}.  Here, we consider a problem where, after a correction from a coarse level in the hierarchy, there are two groups of nodes that have negative approximate solutions.  As a first step in correcting this behaviour, these groups are identified, along with neighboring nodes that have positive values.  In this example, we note that the group of nodes $\{4,5,6\}$ has positive neighboring nodes 3 and 7, and that the group of nodes $\{11,12,13\}$ has positive neighboring nodes 10 and 14.  Assuming the values at points 3, 7, 10, and 14, are more accurate approximations of the true solution, we could now simply solve localized problems to determine values for the points in each group, treating their positive neighbors as Dirichlet boundary conditions, or use a simple interpolation scheme, such as linearly interpolating between the positive neighboring values, as is shown in the figure. To identify such groups of non-positive entries in the solution, we run a simple iteration through all the nodes, and place successive nodes with non-positive solution values into one group.

To perform local corrections using interpolation, we suppose an identified group contains nodes $\mathcal M = \{x_{j}, x_{j+1}, \ldots, x_{j+k} \}$, where $j > 0$ and $j+k < n$.  Then, we correct the values $u_{j+l}$ at nodes $x_{j+l}$ for $0 \leq l \leq k$ locally by linear interpolation from values $u_{j-1}$ and $u_{j+k+1}$, as
\begin{equation}\label{eq:lin-interp-corrrection}
u_{j+l} \leftarrow u_{j-1} + \frac{u_{j+k+1}-u_{j-1}}{x_{j+k+1}-x_{j-1}}(x_{j+l}-x_{j-1}), \quad 0 \leq l \leq k.
\end{equation}
Though it may be possible to define higher-degree polynomial interpolation, or even more complicated positivity-preserving interpolation schemes according to the concept of algebraic smoothness, see e.g. \Cref{sec:AMG}, we find that a simple linear interpolation works well in the various examples that we have considered, so use only this strategy here.

    \definecolor{wrwrwr}{rgb}{0.38,0.38,0.38}
    \begin{figure}
        \centering
        \begin{tikzpicture}[line cap=round,line join=round,>=triangle 45,x=1cm,y=1cm]
    \begin{axis}[
    x=0.6cm,y=0.6cm,
    axis lines=middle,
    xlabel=$x$,
    xmin=-1,
    xmax=18,
    xtick={0,1,2,3,4,5,6,7,8,9,10,11,12,13,14,15,16},
    xticklabels={,,,3,,,,7,,,10,,,,14,,},
    ymin=-2,
    ymax=3,
    ytick=\empty,
    hide y axis]
    
    \addplot+[dashed, color=wrwrwr, mark=*, mark options={solid,fill=green}] 
    coordinates {
        (3,1.2)
        (4,1.3) (5,1.4) (6,1.5) 
        (7,1.6)
    };
    
    \addplot+[dashed, color=wrwrwr, mark=*, mark options={solid,fill=green}] 
    coordinates {
        (10,1.2) 
        (11,1.1) (12,1) (13,0.9) 
        (14,0.8)
    };
    
    \addplot+[color=cyan, mark options={solid,fill=black}] 
    coordinates {
        (0,0) (1,0.5) (2,0.9) (3,1.2)
        (4,-0.8) (5,-1) (6,-0.7) (7,1.6)
        (8,1) (9,1.5) (10,1.2) (11,-0.8)
        (12,-0.2) (13,-1.2) (14,0.8) (15,0.8)
        (16,0)
    };
    
    \draw (0,-0.1) node[anchor=north] {0};
    
    \draw[dashed] (3,0) -- (3,1.2);
    \draw[dashed] (4,0) -- (4,1.3);
    \draw[dashed] (5,0) -- (5,1.4);
    \draw[dashed] (6,0) -- (6,1.5);
    \draw[dashed] (7,0) -- (7,1.6);
    
    \draw[dashed] (10,0) -- (10,1.1);
    \draw[dashed] (11,0) -- (11,1);
    \draw[dashed] (12,0) -- (12,0.9);
    \draw[dashed] (13,0) -- (13,0.8);
    \draw[dashed] (14,0) -- (14,0.8);
    
    \end{axis}
    \end{tikzpicture}
        \caption{Local corrections on an example approximate solution}
        \label{fig:unigrid_local_correction}
    \end{figure}
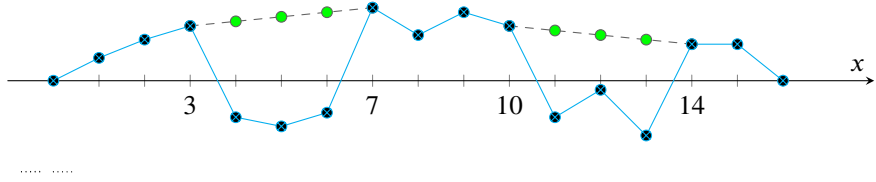

    Using simple linear interpolation to correct solution values is feasible for one-dimensional problems, but does not generalize effectively to higher dimensions.  As an alternative approach, we consider correcting negative solution values using local relaxation to restore positivity.  In this approach, we still identify those points where the approximate solution after a unigrid correction step is negative, but now put them all into one group only.  To correct the negative values on the problematic points, we fix the positive corrected values in the approximation, and re-solve the subset of solution with problematic (negative) values through a simple iteration over this block.  In practice, we find that applying point-wise Gauss-Seidel iterations repeatedly (lexicographically within the group of negative points) is quite effective.  That is, given a set $\mathcal M$ with the indices of the problematic points, local fine-grid Gauss-Seidel relaxation on those points,
\begin{equation}\label{eq:GS-correction}
u_i \leftarrow (b_i - \sum_{j\neq i}A_{ij}u_j)/A_{ii}, \quad \forall i\in\mathcal M,
\end{equation}
is repeated until the solution approximation at all points becomes positive again.  In practice, we reexamine $\mathbf u$ after each relaxation sweep in~\cref{eq:GS-correction}, to update the definition of $\mathcal M$, and only relax on points where $u_i < 0$ on subsequent sweeps.

The unigrid algorithm with local correction is
given in \cref{alg:unigrid_local}. We note that this has the potential to become
costly if large groups of points become negative in one correction step.  In our
tests, discussed below, we found that this cost did not become substantial.

\begin{algorithm}
  \begin{algorithmic}[1]
    \Procedure{UG\_SOLVE\_CORRECTION}{$A, \mathbf b, \mathbf u, \nu, \{I_k^0\}$}
\For{$k = 0, 1, ..., l$} 
  \For{$i = 1, ..., \nu$} 
     \For{$j = 1, ..., n_k$} 
     \State Compute correction ${\delta}^{(k)}_j = \frac{\langle \mathbf b-A\mathbf u, I_k^0\mathbf{e}_j^{(k)} \rangle}{\langle AI_k^0\mathbf{e}_j^{(k)}, I_k^0\mathbf{e}_j^{(k)}\rangle}$
      \State Update approximation $\mathbf u = \mathbf u + {\delta}^{(k)}_jI_k^0\mathbf{e}_j^{(k)}$
       \If{$\min\{ \mathbf u\} \leq 0$ }
           \State Iterate through $\mathbf u$ and find all non-positive groups $\mathcal M = \{\mathcal{M}_1, \mathcal{M}_2, \ldots, \mathcal{M}_s \}$
           \If{Linear interpolation}
               \For{each identified group $\mathcal{M}_j\in\mathcal M$}
               \State Perform local correction by linear interpolation, as in~\cref{eq:lin-interp-corrrection} 
               \EndFor
           \EndIf
           \If{Gauss-Seidel correction}
               \While{$\min\{u_j\}_{j\in\mathcal M} \leq 0$}
               \State Perform local correction by Gauss-Seidel relaxation, as in~\cref{eq:GS-correction}
               \State Update $\mathcal M$
               \EndWhile
           \EndIf
      \EndIf 
   \EndFor 
 \EndFor 
\EndFor
\State Return $\mathbf u$
\EndProcedure
\caption{Unigrid Method with Local Correction}\label{alg:unigrid_local}
\end{algorithmic}
\end{algorithm}

\section{Numerical Experiments}\label{sec:numerical_experiments}

We now consider several one- and two-dimensional model problems to demonstrate both the performance of the algorithms discussed above and the potential for non-positive approximations to be generated by standard multigrid algorithms, such as Ruge-St\"uben AMG.  All examples are coded in Matlab, and the same code is used to generate the AMG interpolation and coarse-grid operators for both the direct AMG algorithm and the unigrid methods.  While AMG can also be used as a preconditioner, we note that the use of (flexible) Krylov acceleration would negate the positivity-preserving properties of the unigrid algorithms, so we consider only stationary iterations.  Because the codes are written in Matlab and not optimized, we do not report direct algorithmic timings but, instead, note the additional costs of enforcing positivity through appropriate work measures on the unigrid schemes.

\subsection{1D Linear Problem}\label{sec:simple_1d_ex}
We first consider the simple model problem
\[
-(\sigma u')' = \sin(\pi x),
\]
on the interval $0 < x < 1$ with boundary conditions $u(0) = u(1) = 0$, and diffusion coefficient
\begin{equation*}
    \sigma = 
    \begin{cases}
    10^{12}, & 0 < x < 0.4,\\
    1,    & 0.4 \leq x < 1.\\
    \end{cases}
\end{equation*}
We note that this example is constructed specifically to be difficult to preserve positivity of the approximate solution by classical methods, because coarse-grid correction will naturally couple values across the jump in the diffusion coefficient that is (purposefully) not aligned with traditional coarse-grid points for geometric multigrid.

\begin{figure}[h]
    \centering
    \includegraphics[width=0.47\textwidth]{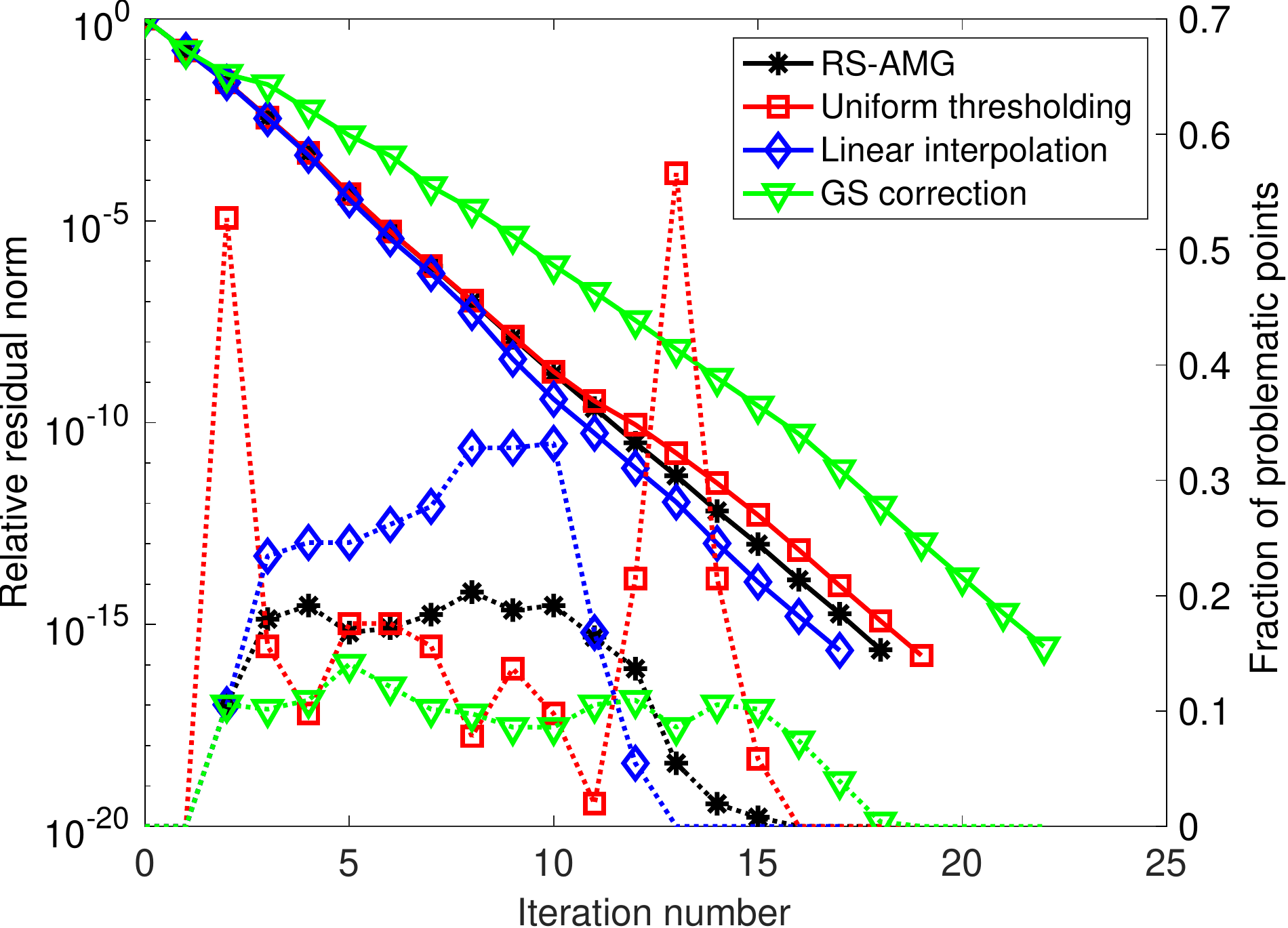}
    \includegraphics[width=0.47\textwidth]{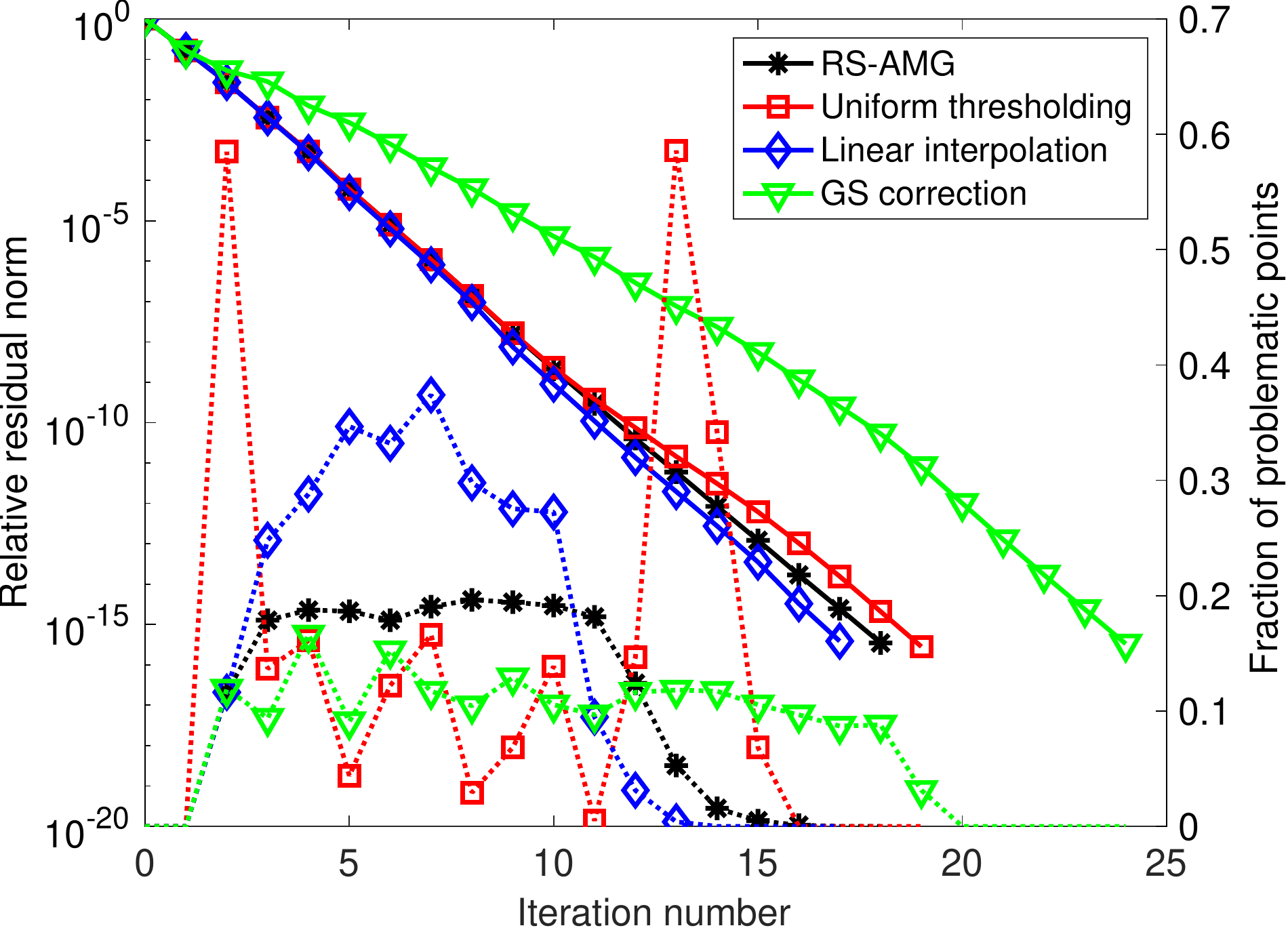}
    \caption{Numerical results for 1D linear model problem with $N=256$ (at left) and $N=1024$ (at right).  Solid lines depict residual convergence histories (left axis) while dotted lines (right axis) depict the number of points (as a fraction of the finest-grid size) where either negative entries were present after the iteration (for RS-AMG) or where a unigrid correction was adapted to preserve positivity (for all unigrid results).}
    \label{fig:1dlinear}
\end{figure}

\Cref{fig:1dlinear} presents convergence data on meshes with $N$ elements for $N=256$ (at left) and $N=1024$ (at right), with a stopping tolerance of a relative reduction in the Euclidean norm of the residual by a factor of $10^{15}$.  As a baseline, we consider the performance of standard RS-AMG, which converges in 18 iterations for both grid sizes.  The corresponding dotted line shows the fraction of the total mesh points that have a negative value after each iteration.  For both grid sizes, we see that about 20\% of the grid points have such negative values over a range of about 10 of these iterations, and there remain negative values in the solution vector until the solution is nearly converged.  (We note also that the $10^{15}$ relative reduction in the residual is quite a strict stopping tolerance; using a more typical relative residual reduction tolerance of $10^8$ or so would yield RS-AMG convergence to a solution with negative entries.)

In comparison, we consider three variants on the unigrid algorithm:
\begin{itemize}
\item Using uniform thresholding, as in~\cref{alg:unigrid_uniform},
\item Using local correction based on linear interpolation between neighboring positive values, as in~\cref{alg:unigrid_local}, and
\item Using local correction based on Gauss-Seidel relaxation on non-positive values, as in~\cref{alg:unigrid_local}.
\end{itemize}
For each of these, we similarly measure the work used to preserve positivity at each unigrid relaxation step.  For uniform thresholding and local correction using linear interpolation, $\delta_j^{(k)}I_k^0\vec{e}_j^{(k)}$ is computed after each unigrid relaxation step, and we count the number of points that would be made negative by this correction to $\vec{u}$ that are ``recovered'' by either the uniform thresholding (damping the correction to avoid negative values in $\vec{u}$) or by taking the step and then correcting negative values based on linear interpolation from neighboring positive values.  We sum these over all relaxation steps in the unigrid iteration, and express the total number of points recovered in this way normalized by grid size, as a ``fraction of problematic points'', noting that it is possible for this fraction to be greater than 1 if the correction steps occur many times in a cycle.  Similarly, for the local correction based on Gauss-Seidel relaxation, we count the number of times that a Gauss-Seidel correction is applied to a point for each step of relaxation within a unigrid sweep, noting that the same point may be corrected more than once in a unigrid step if its value is not made positive by a single step of Gauss-Seidel relaxation.  Again, this is reported in the figures, normalized by problem size, as the ``fraction of problematic points'' for the algorithm.  In all plots, we use solid lines to show convergence of the algorithm (reduction in the Euclidean norm of the residual relative to its initial value, shown on the left-hand $y$-axis), and dashed lines to report this fraction of problematic points (shown on the right-hand $y$-axis).

The unigrid data in~\cref{fig:1dlinear} shows that using the unigrid methodology has only a small effect on the multigrid convergence rates, with uniform thresholding and local correction using linear interpolation both tracking the RS-AMG convergence curve, and reaching the (strict) convergence tolerance within one iteration of RS-AMG.  Convergence of the method using GS correction is somewhat slower, requiring 22 iterations for $N=256$ and 24 for $N=1024$.  We suspect this is caused by the effects of local Gauss-Seidel relaxation on large-scale unigrid corrections (from coarse levels in the hierarchy), where algorithmic convergence is being sacrificed in favor of non-negativity of the evolving approximate solution.  Considering the work done by the unigrid algorithm to ensure positivity, we note that none of the algorithms do substantial work in comparison to the cost of a fine-grid Gauss-Seidel sweep.  Using uniform thresholding, we find that generally a small percentage of the unigrid relaxation steps are impacted (totaling up to 20\% of the points in the grid for all but 3 iterations), but that substantial damping is used on some iterations - this is likely from damping being applied to very coarse-scale corrections.  For the local correction approaches, we see that linear interpolation generally corrects a larger fraction of points than GS correction, but that neither has substantial cost, particularly when averaged over many iterations.  For the GS correction, for example, the total cost of all correction steps is just under that of 2 full Gauss-Seidel sweeps on the finest grid.

\subsection{Nonlinear grid generation in 1D}
We next consider a nonlinear grid generation problem in 1D, following Huang and Russell~\cite{WHuang_RDRussell_2011a}.  For this problem, we consider the nonlinear 1D diffusion equation,
\begin{equation}
    -\frac{d}{dx}\bigg( a(u(x)) \frac{du}{dx} \bigg) = 0, \quad x\in (0,1),
\end{equation}
with $u(0)=0$ and $u(1)=1$, where the function $u(x)$ determines a computational grid on the interval $[0,1]$ with $N$ elements by sampling $u(x)$ at the $N+1$ evenly-spaced points $u(j/N)$ for $0\leq j \leq N$.  Here, the coefficient $a(u)$ determines the mesh density; we choose the example
\begin{equation}
a(u) = \begin{cases} 1000 & \text{if }u < 0.5 \\ 1 & \text{if }u>0.5 \end{cases},
\end{equation}
which leads to greater mesh density for small values of $u(x)$.

Here, we follow the simplest form of the Picard iteration from Section~\ref{ssec:discretization}, solving
\begin{equation}\label{eq:nldiff_picard}
    A(\mathbf u^{(k-1)})\mathbf u^{(k)} = \mathbf b^{(k-1)},
\end{equation}
at each iteration, $k$, where $A(\vec{u}^{(k-1)})$ is the $(N-1)\times (N-1)$ centered finite-difference discretization matrix over the interior points of the mesh for fixed (previous) approximation $\vec{u}^{(k-1)}$, and $\vec{b}^{(k-1)}$ is a non-negative vector determined by the fixed boundary data for $u(0)$ and $u(1)$ (but indirectly dependent on the evolving entries in $\vec{u}$ from the ``eliminated'' off-diagonal entries in rows adjacent to the boundary points).  Since $A(\vec{u}^{(k-1)})$ is an irreducible M-matrix and $\vec{b}^{(k-1)}$ has one positive entry (since $u(0)=0$), we know that $\vec{u}^{(k)}$ determined by~\cref{eq:nldiff_picard} should be positive at each iteration.  We use an initial guess of $\vec{u}^{(0)} = \vec{x}$, the vector whose $i^\text{th}$ component is $x_i = i/N$, corresponding to a uniform grid function, $u(x)$.  We use a relative tolerance on the nonlinear residual as a stopping criterion on the outer (Picard) loop, requiring that
\begin{equation}\label{eq:nonlinear_stopping_example}
\|\vec{b}^{(k)} - A(\vec{u}^{(k)})\vec{u}^{(k)}\| \leq 10^{-10}\|\vec{b}^{(0)} - A(\vec{u}^{(0)})\vec{u}^{(0)}\|.
\end{equation}
For the inner (AMG or unigrid) iterations, we require that the Euclidean norm of the linear residual be reduced either by a relative factor of $10^8$ (with initial guess for $\vec{u}^{(k)}$ being equal to the final value of $\vec{u}^{(k-1)}$) or below one tenth of the value of the nonlinear stopping criterion on the right-hand side of~\cref{eq:nonlinear_stopping_example}.

\begin{figure}[h]
    \centering
    \includegraphics[width=0.47\textwidth]{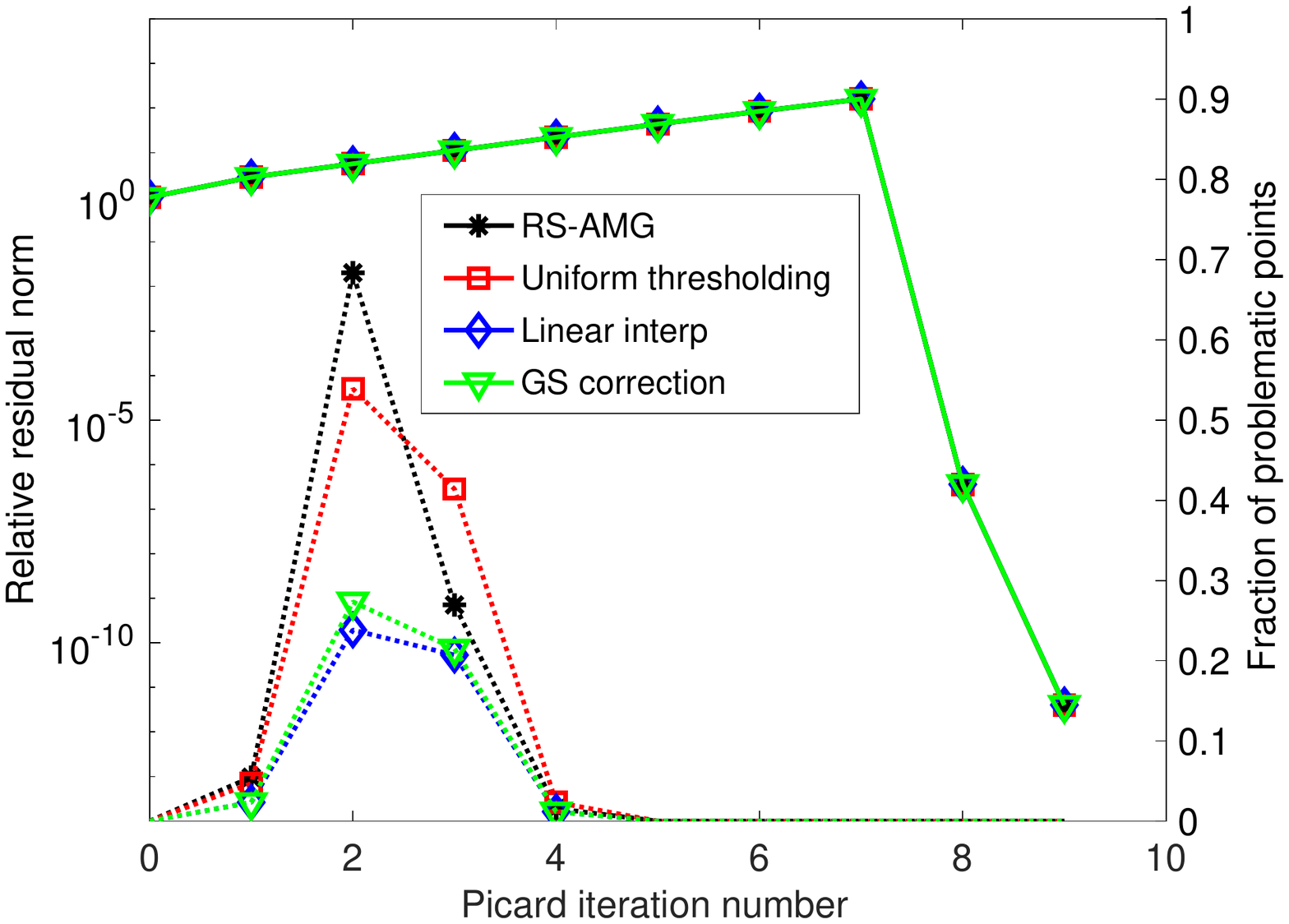}
    \includegraphics[width=0.47\textwidth]{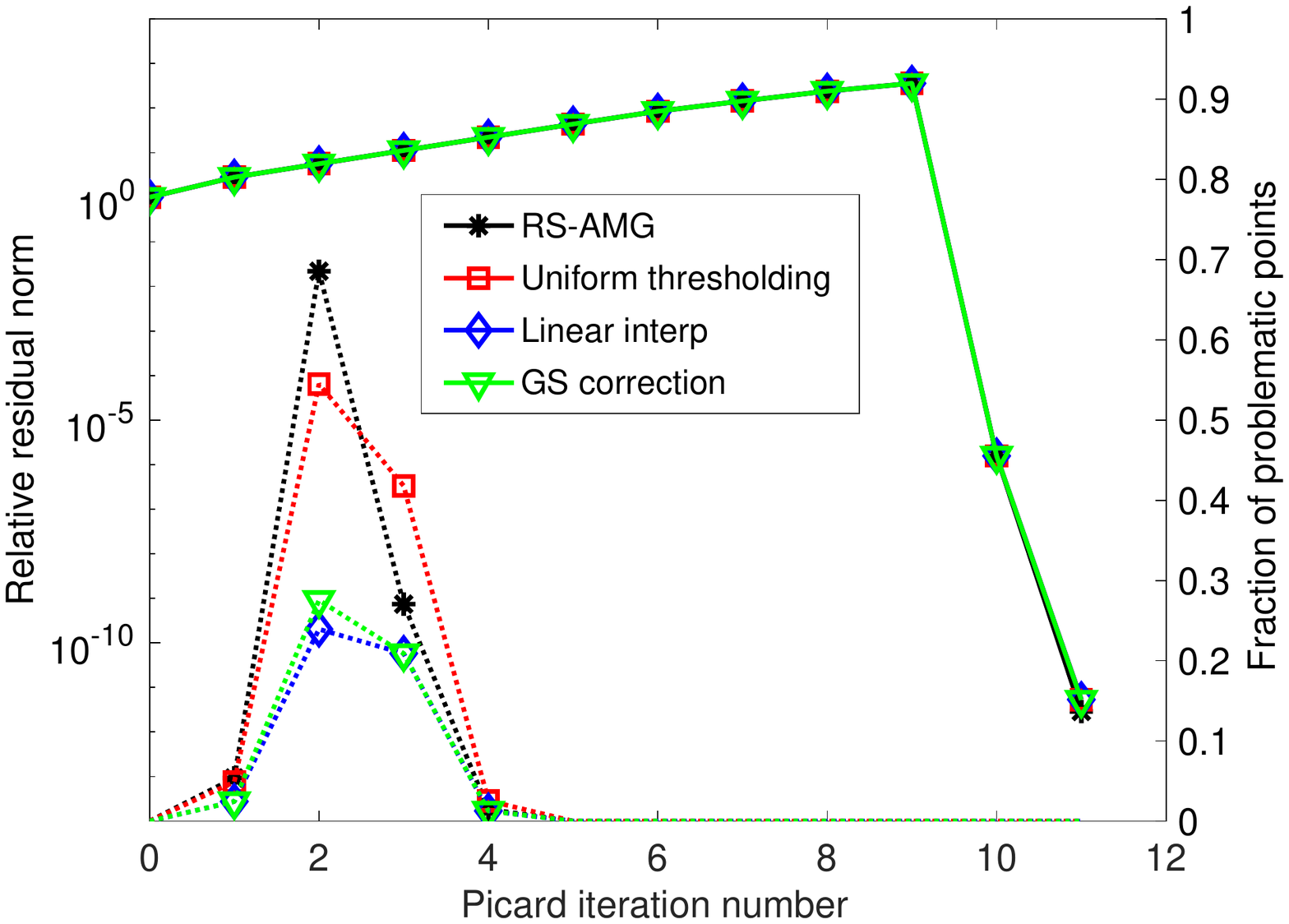}
    \caption{Numerical results for 1D nonlinear grid generation model problem with $N=256$ (at left) and $N=1024$ (at right).  Solid lines depict nonlinear residual convergence histories (left axis) while dotted lines (right axis) depict the number of points (as a fraction of the finest-grid size) where either negative entries were present after a linear iteration (for RS-AMG) or where a unigrid correction was adapted to preserve positivity (for all unigrid results), accumulated over all linear solves within each Picard iteration.}
    \label{fig:1dnonlinear}
\end{figure}

\Cref{fig:1dnonlinear} shows the Picard iteration convergence history for this example.  We note that all methods show very similar nonlinear convergence, with a steady growth in the residual norm through a number of iterations, followed by effective convergence in the final iterations.  This is due to the nature of the solution of the nonlinear problem, which has a boundary-layer like solution, with near-linear growth from $u=0$ at the left-hand endpoint to $u\approx 0.2$ near the right-hand endpoint, then quick growth to satisfy the right-hand boundary condition $u(1)=1$.  The Picard iterates slowly adapt to this, ``seeing'' more gridpoints in the $u<0.5$ region with each iteration until the diffusion coefficient is properly resolved, then the nonlinear solution converges.  As this happens, we see several Picard iterations where there are significant numbers of negative points in the RS-AMG iterates, or where the unigrid algorithms take steps to ensure positivity of the solutions.  Here, in some contrast to the earlier example, we see that RS-AMG yields many negative solution values at the second iterate (primarily in the first two linear iterations).  Among the unigrid methods, uniform correction does the most work to preserve positivity (again in the first two linear iterations per Picard iteration), while linear interpolation and GS correction do substantially less work to preserve positivity.

\begin{figure}[h]
    \centering
    \includegraphics[width=0.47\textwidth]{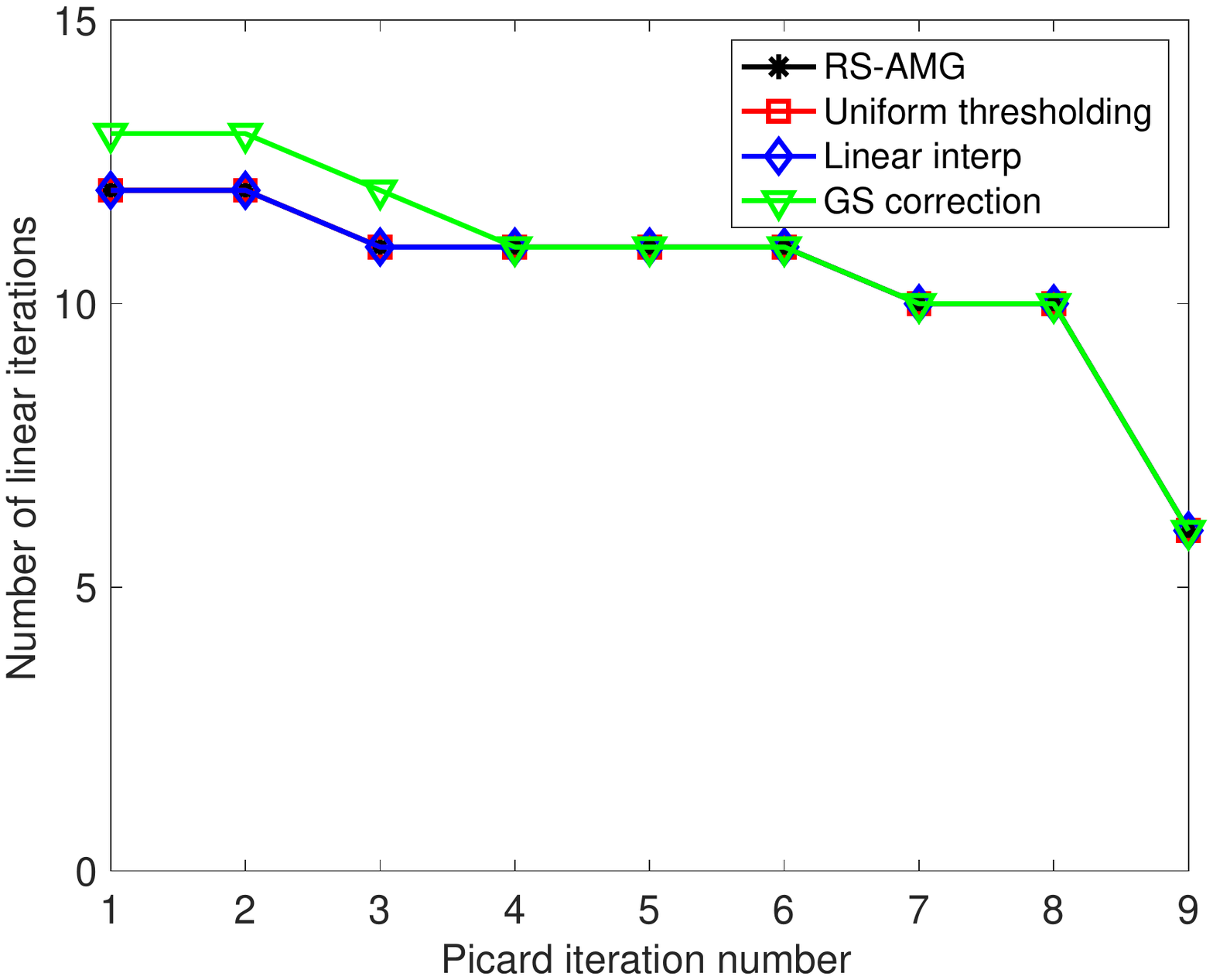}
    \includegraphics[width=0.47\textwidth]{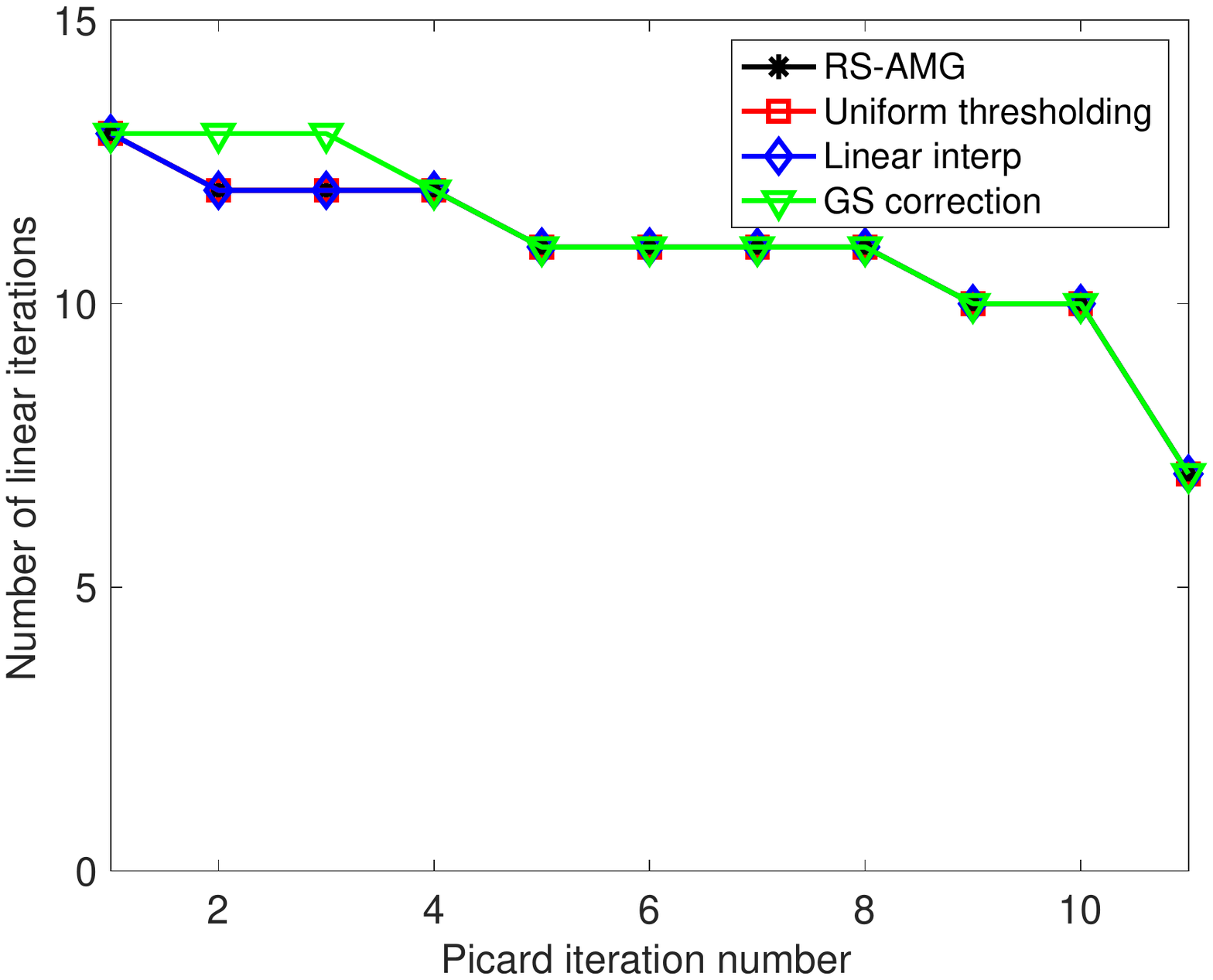}
    \caption{Number of linear iterations for each Picard iteration for the 1D nonlinear grid generation model problem with $N=256$ (at left) and $N=1024$ (at right).}
    \label{fig:1dnonlinear_linear}
\end{figure}

\Cref{fig:1dnonlinear_linear} shows the number of linear iterations required for each method to reach the convergence criterion.  In general, we see there is little variation between the methods, with unigrid using GS correction requiring one additional linear iteration in the early nonlinear solves.  The number of iterations drops substantially in the final iterations, as the linear stopping criterion becomes looser when we are close to convergence.  The allows us to conclude that the additional costs in the unigrid approach are small, in comparison to the overall solve cost of about 10 Picard iterations, each with about 10 linear iterations to solution.

\subsection{2D Piecewise-Constant Diffusion}

We now consider a two-dimensional problem,
\begin{equation*}
    -\nabla \cdot(\sigma \nabla u) = \sin(\pi xy),
\end{equation*}
with
\begin{equation*}
    \sigma = 
    \begin{cases}
    10^6, & 0 < x < 0.8,\; 0 < y < 0.6\\
    1,    & \text{otherwise}.\\
    \end{cases}
\end{equation*}
We use a bilinear finite-element discretization on uniform $N\times N$ meshes of the unit square domain, with homogeneous Dirichlet boundary conditions on all four edges of the domain.  The finite-element stiffness matrix is evaluated exactly, while the right-hand side vector is calculated using midpoint quadrature on each element.  For all methods below, we use the same initial guess for the solution, as a constant vector with values $0.1$.

\begin{figure}[h]
    \centering
    \includegraphics[width=0.47\textwidth]{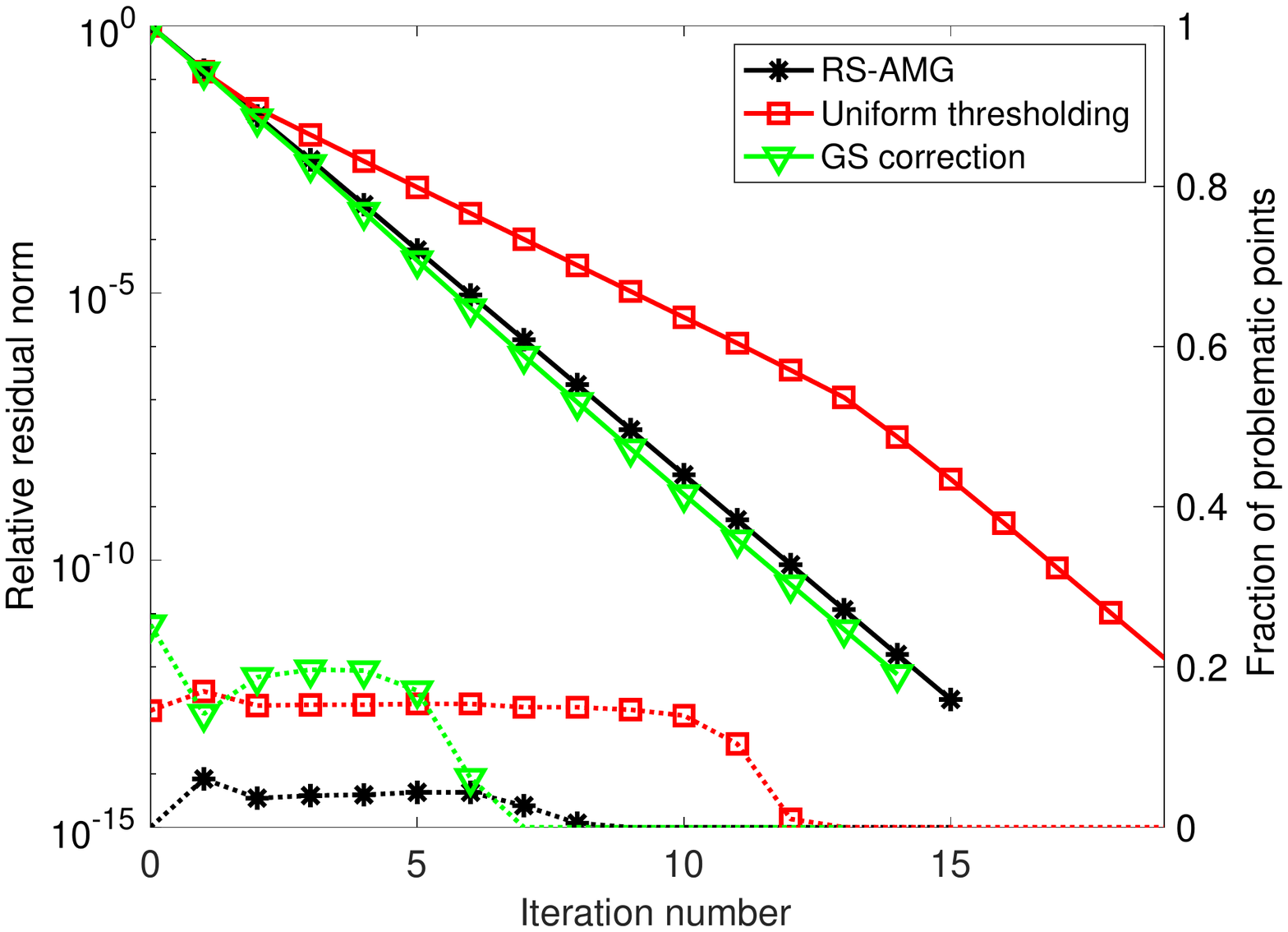}
    \includegraphics[width=0.47\textwidth]{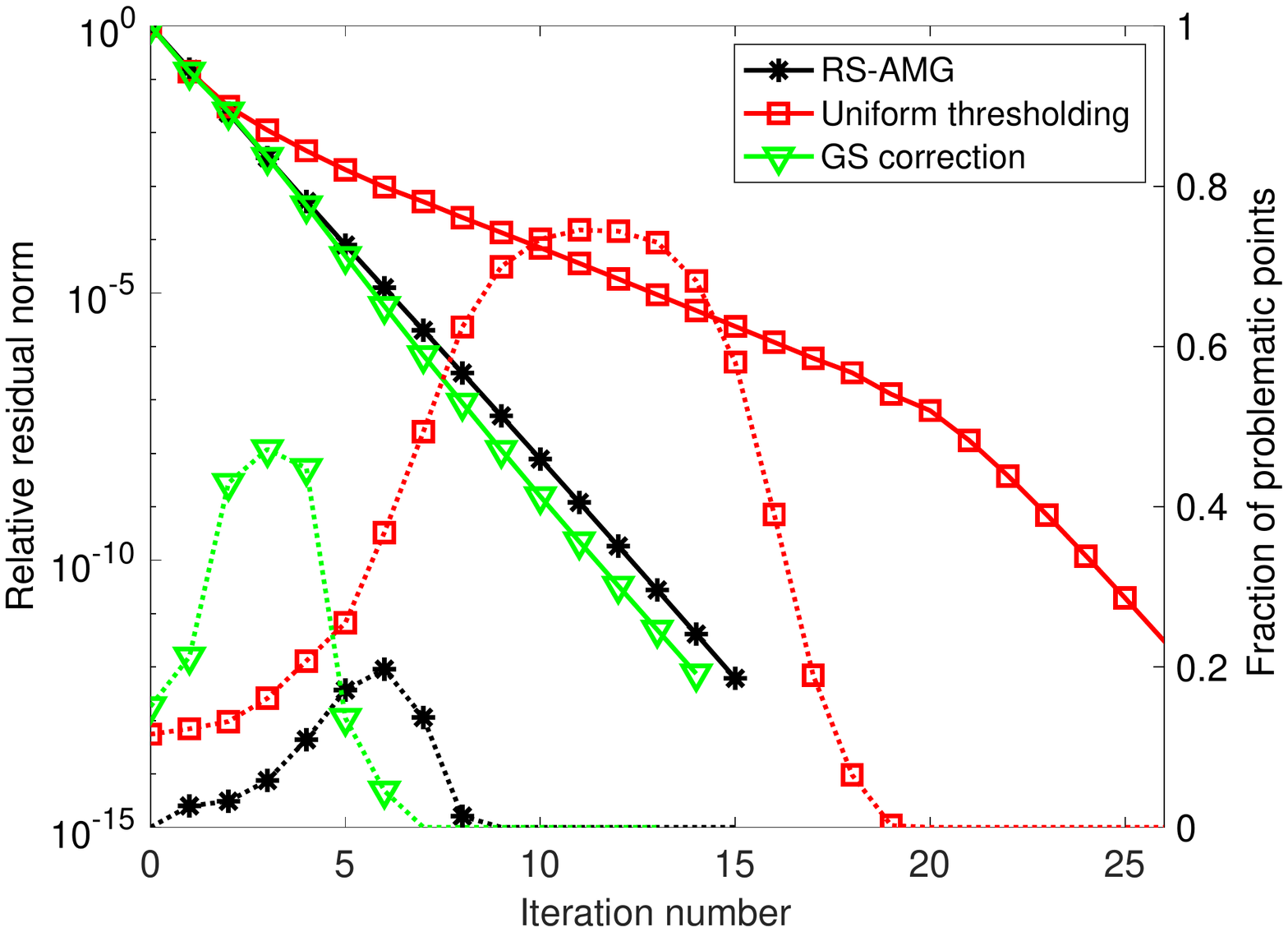}
    \caption{Convergence history for two-dimensional piecewise-constant diffusion problem on $N\times N$ grids  for $N=32$ (left) and $N=64$ (right).  Solid lines denote residual convergence histories (left axis) while dotted lines (right axis) depict the number of points (as a fraction of the finest-grid size) where either negative entries were present after the iteration (for RS-AMG) or where a unigrid correction was adapted to preserve positivity (for all unigrid results).}
    \label{fig:2d_linear_first}
\end{figure}

\Cref{fig:2d_linear_first} depicts convergence histories for this problem for $N=32$ and $64$.  We note that the unigrid method with local correction using linear interpolation does not easily generalize to two-dimensional problems, so we no longer include it in the comparison.  We also note that the RS-AMG results are more-or-less as expected, showing convergence in 15 iterations for both grid sizes, with some non-positive entries arising in the early iterations.  There is a noticeable, but small, increase in the number of non-positive entries produced by RS-AMG for $N=64$ in comparison to $N=32$.  The use of uniform thresholding with unigrid, in contrast, shows both an increasing number of iterations (19 for $N=32$ and $26$ for $N=64$) and increasing damping used as we increase $N$.  Thus, while the asymptotic convergence of unigrid with uniform thresholding is clearly similar to that of RS-AMG, it is a much less efficient solver.  In contrast, unigrid with local GS correction appears to converge almost identically to RS-AMG, in 14 iterations for both grid sizes.  There is again some noticeable growth in the total number of local GS steps performed relative to grid size as $N$ increases, but it still comes at a total cost less than that of a handful of sweeps of Gauss-Seidel on the entire finest-grid problem.

\subsection{Checkerboard diffusion equation in 2D}

\begin{figure}
  \centering
  \includegraphics{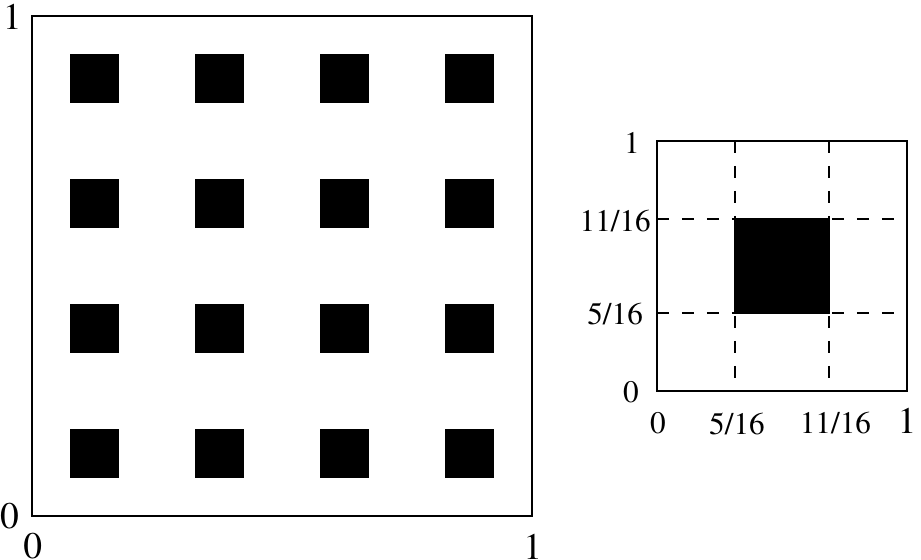}
  \caption{At left, sketch of the diffusion coefficient given in~\cref{eq:checkerboard} for the case of $p=4$, with dark regions showing where the diffusion coefficient is 1, and light regions denoting where the diffusion coefficient is 1000.  At right, a sketch of the unit cell that is scaled and tiled into a $p\times p$ array to determine $\sigma$.}\label{fig:checkerboard}
\end{figure}

For our final example, we again consider the bilinear finite-element discretization of the two-dimensional problem,
\begin{equation*}
    -\nabla \cdot(\sigma \nabla u) = \sin(\pi xy),
\end{equation*}
on uniform $N\times N$ meshes of the unit square with homogeneous Dirichlet boundary conditions.
We now consider a so-called ``checkerboard'' diffusion coefficient, with periodicity $p$, taking
\begin{equation}\label{eq:checkerboard}
    \sigma = 
    \begin{cases}
    1, & 5/16 < [px] < 11/16 \text{ and } 5/16 < [py] < 11/16\\
    1000,    & \text{otherwise}\\
    \end{cases},
\end{equation}
where we use $[px]$ to denote the fractional part of real number $px$.  In the numerical examples below, we take $p = N/16$, so that this prescribes a piecewise constant diffusion coefficient on the finest grid, with jumps in the diffusion coefficient value aligned with the mesh cells.  A sketch of the coefficient for the case of $N=64$ and $p=4$ is shown in~\cref{fig:checkerboard}; for larger values of $N$ and $p$, the pattern repeats more frequently.

\begin{figure}[h]
    \centering
    \includegraphics[width=0.47\textwidth]{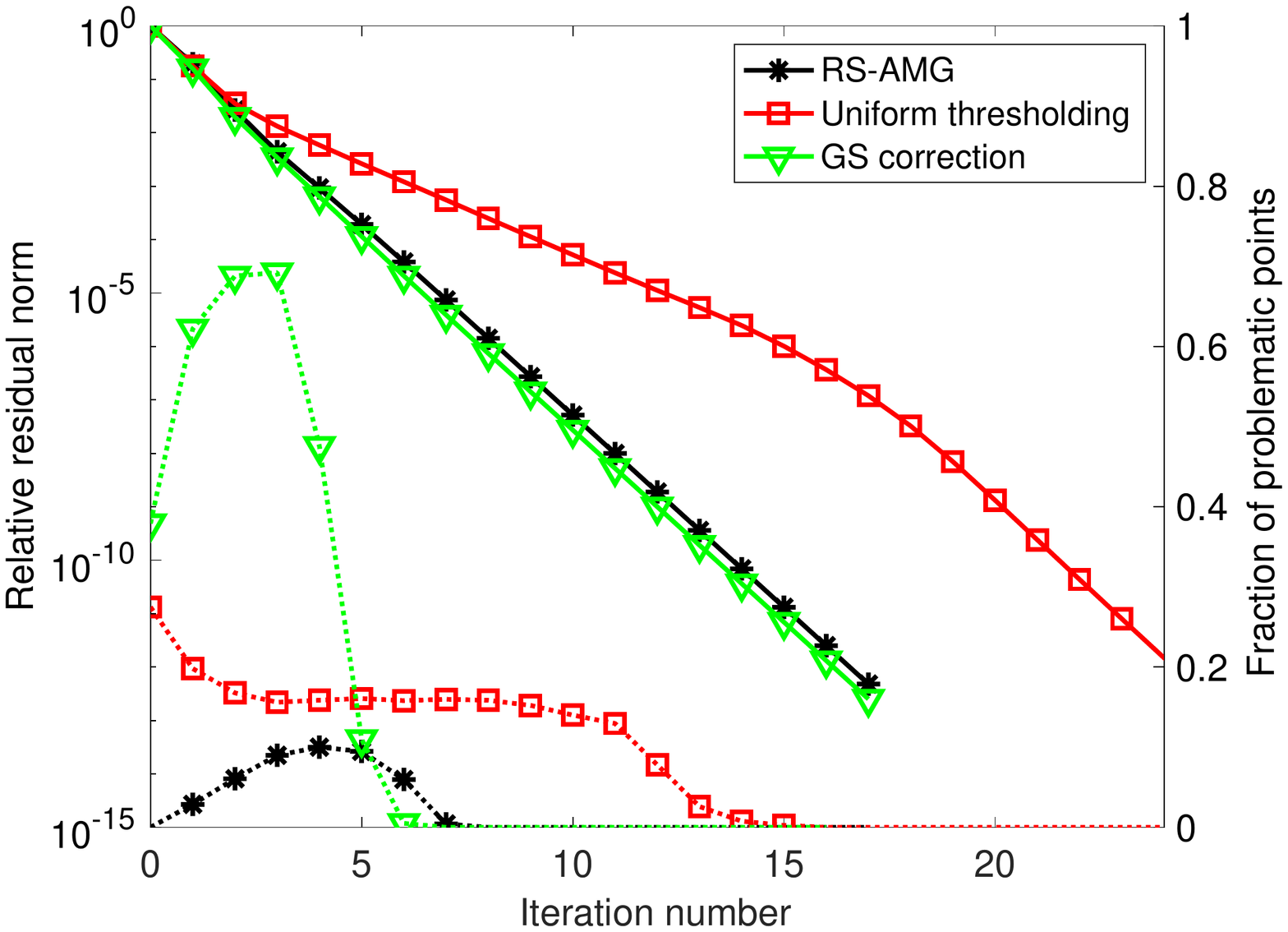}
    \includegraphics[width=0.47\textwidth]{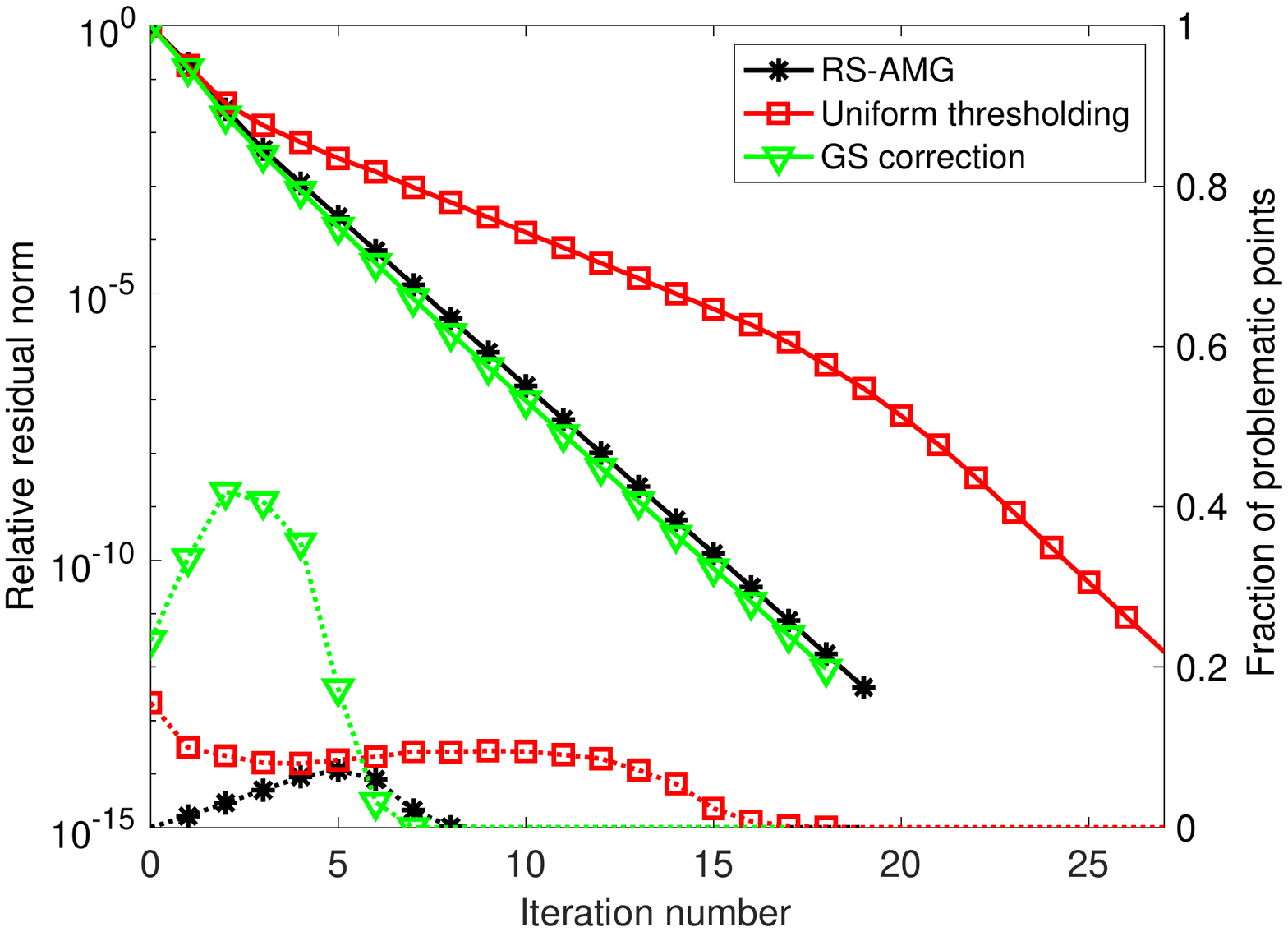}
    \caption{Convergence history for two-dimensional checkerboard diffusion problem on $N\times N$ grids  for $N=128$ (left) and $N=256$ (right).  Solid lines denote residual convergence histories (left axis) while dotted lines (right axis) depict the number of points (as a fraction of the finest-grid size) where either negative entries were present after the iteration (for RS-AMG) or where a unigrid correction was adapted to preserve positivity (for all unigrid results).}
    \label{fig:2d_checkerboard}
\end{figure}

\Cref{fig:2d_checkerboard} depicts convergence results for this problem with $N=128$ ($p=8$) and $N=256$ ($p=16$), starting from an initial guess for the solution vector of all ones.  Broadly, these results look similar to those in~\cref{fig:2d_linear_first}, with slower convergence of the unigrid algorithm using uniform thresholding than GS correction.  Notably, both unigrid with GS correction and RS-AMG encounter substantial numbers of negative entries in the approximation vector over the first 7 or so iterations, while unigrid using uniform thresholding applies thresholding over nearly double as many iterations.  This, in turn, leads to a visibly slower rate of residual convergence of the algorithm using uniform thresholding.  In contrast, both RS-AMG and unigrid with GS correction converge in roughly the same number of iterations, with only a slight degradation as the problem size increases.  Again, the additional cost of the GS correction steps is bounded by that of a few fine-grid Gauss-Seidel iterations, indicating there is little practical cost of maintaining positivity of the approximate solution.

\section{Conclusions}\label{sec:conclusion}

The solutions of many PDEs are known to satisfy additional constraints on their properties that are often lost in the numerical discretization and solution processes.  Here, we propose and test several variants of unigrid that are applied to problems with natural pointwise positivity constraints that are preserved through discretization by M-matrix structure.  The two most general variants are using thresholding to damp corrections that would introduce negative entries to the approximate solution and using localized Gauss-Seidel iterations to correct such negative entries after they are introduced.  Numerical results for two-dimensional variable-coefficient diffusion problems show that the thresholding approach leads to increased cost and iteration counts in comparison to a standard solver, like RS-AMG, whose performance is closely matched by the GS-correction variant.

Future work includes developing an efficient, parallel implementation in a
compiled language to support testing on a wider range of problems, including
three-dimensional PDEs. The main challenge to be addressed in such work is to
efficiently identify problematic points after each unigrid correction, and to
synchronize such corrections across MPI ranks.  An exploration of time dependent
applications such as thin-film problems \cite{bertozzi1998,Zhornitskaya1999-lj,matta2011} is also a natural next step.

\bibliography{unigrid}

\section*{Acknowledgments}
This work was partially supported by NSERC Discovery Grants to RH and SM.

\subsection*{Conflict of interest}

The authors declare no potential conflict of interests.

\end{document}